\documentclass[journal,twoside,web]{ieeecolor}

\usepackage{generic}  % Change this to TAC (find correct code)
\usepackage{cite}
\usepackage{amsmath,amssymb,amsfonts}
\usepackage{algorithmic}
\usepackage{graphicx}
\usepackage{textcomp}
\usepackage{mathtools}
\usepackage{hyperref}

\newtheorem{Theorem}{Theorem}
\newtheorem{Lemma}{Lemma}
\newtheorem{Problem}{Problem}
\newtheorem{Remark}{Remark}

\newtheorem{Assumption}{Assumption}
\newtheorem{Definition}{Definition}

\DeclareMathOperator{\R}{\mathbb R}

\DeclareMathOperator*{\argmin}{arg\,min}

\newcommand{\bb}{\boldsymbol}

\def\BibTeX{{\rm B\kern-.05em{\sc i\kern-.025em b}\kern-.08em
    T\kern-.1667em\lower.7ex\hbox{E}\kern-.125emX}}
\begin{document}

\title{Consolidated Control Barrier Functions: Synthesis and Online Verification via Adaptation under Input Constraints}

% Authors, Affiliations, and Acknowledgements
\author{Mitchell Black, \IEEEmembership{Student Member, IEEE}, and Dimitra Panagou, \IEEEmembership{Senior Member, IEEE}
\thanks{This paper was submitted for review on March 7, 2023. The authors would like to acknowledge the support of the National Science Foundation award number 1931982.}
\thanks{Mitchell Black is with the Department of Aerospace Engineering, University of Michigan, Ann Arbor, MI 48109 USA (e-mail: mblackjr@umich.edu). }
\thanks{Dimitra Panagou is with the Department of Robotics, University of Michigan, Ann Arbor, MI 48109 USA, and the Department of Aerospace Engineering (e-mail: dpanagou@umich.edu).}
}

\maketitle

\begin{abstract}
In this paper, we develop a novel adaptation-based approach to constrained control design under multiple state and input constraints. Specifically, we introduce a method for synthesizing any number of time-varying candidate control barrier functions (CBF) into one consolidated CBF (C-CBF) candidate, and propose a predictor-corrector optimization-based adaptation law for the weights of the constituent constraint functions that certifies the C-CBF as valid for a class of nonlinear, control-affine systems. We prove this result by showing that the adapted weights are guaranteed to confer sufficient control authority to meet the new, adaptive C-CBF condition in perpetuity despite input constraints, which thereby permits its use in a quadratic program based control law. We then illustrate the performance of our controller on an academic example, and further highlight that it is successful even for constraint functions with higher or mixed relative-degree by simulating a reach-avoid problem for bicycle robots, which we use to demonstrate how our approach out-performs two baseline approaches.
\end{abstract}


\section{Introduction}\label{sec.introduction}

Since the arrival of control barrier functions (CBFs) to the field of safety-critical systems \cite{wielandallgower2007cbf}, much attention has been devoted to the development of their viability for safe control design \cite{ames2017control,Xiao2019HOCBF,cortez2019cbfmechsys}. As a set-theoretic approach founded on the notion of forward invariance, CBFs certify adherence to constraints in that they ensure that any state beginning within a given set remains so for all future time. In the context of control design, CBF conditions are often used as constraints in quadratic program (QP)-based control laws, either as safety filters \cite{Chen2018Obstacle} or in conjunction with stability or liveness constraints (e.g., control Lyapunov functions) \cite{Garg2021Robust}. Their utility has been successfully demonstrated for a variety of safety-critical applications, including mobile robots \cite{Chen2021Guaranteed,Jankovic2021Collision}, unmanned aerial vehicles (UAVs) \cite{Xu2018Safe,Khan2020Cascaded}, and autonomous driving \cite{Black2022ffcbf,Yaghoubi2021RiskBoundedCBF}. But while it is now well-established that synthesizing a CBF for a constraint set serves as a certificate of constraint adherence via set invariance, the verification of \textit{candidate} CBFs as \textit{valid} is in general a challenging problem.

\begin{figure*}
    \centering
        \includegraphics[trim=0.4cm 21.8cm 0.1cm 0cm,clip,width=1\textwidth]{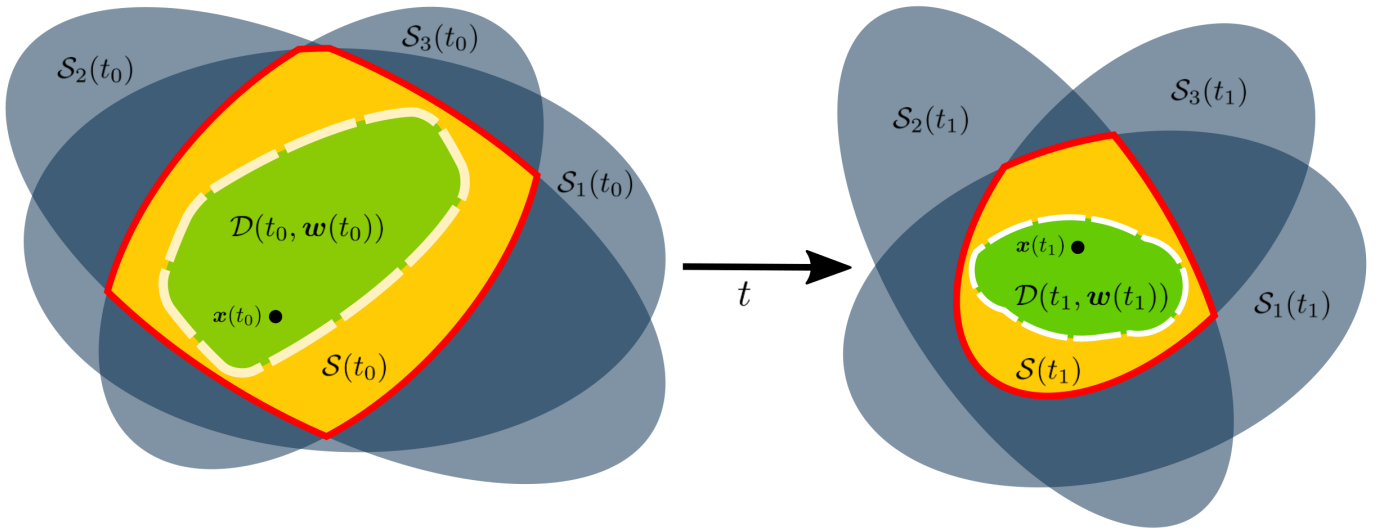}
    \caption{Parameter adaptation for our consolidated control barrier function generates a time-varying, weight-dependent controlled-invariant set $\mathcal{D}(t, \bb{w}(t)) \subset \mathcal{S}(t) = \bigcap_{i=1}^c \mathcal{S}_i(t)$ for a collection of $c$ time-varying constraint sets $\mathcal{S}_i(t)$. The sets shown here evolve across time from the configuration on the left at $t_0$ to that on the right at $t_1$, but at any given time $t$ the state $\bb{x}(t)$ must lie within $\mathcal{D}(t, \bb{w}(t))$.}\label{fig.adaptive_safe_set}
    \vspace{-3mm}
\end{figure*}

% \begin{figure*}
%     \centering
%     \includegraphics[width=0.8\textwidth]{img1.jpg}
%     \caption{This is a wide figure that spans both columns.}
%     \label{fig:img1}
% \end{figure*}

Though for a single candidate CBF there exist guarantees of validity under certain conditions for systems with unbounded \cite{ames2017control} control authority, verifying a candidate CBF under input constraints poses significant challenges. In response, various works have demonstrated the success of verification tools (e.g., sum-of-squares optimization \cite{Wang2018Permissive,clark2021verification}, linear programming \cite{pond2022fast}) in synthesizing a valid CBF offline prior to deployment. These verification certificates, however, may be invalidated by unmodelled phenomena like exogenous disturbances or environmental changes. To address this drawback, the authors of \cite{Breeden2021InputConstraints} propose an online method for guaranteed constraint satisfaction specific to high-order (HO-) CBFs, though their approach typically requires forward simulation of the system trajectories. Alternatively, in \cite{Xiao2022Adaptive} an adaptation-based approach is introduced for guaranteed feasibility of a HO-CBF-QP control law provided that the parameters adapt sufficiently quickly. The above results, however, break down in the presence of multiple constraints.

% While the authors of \cite{Breeden2021InputConstraints} propose a solution to this problem specific to high-order (HO-) CBFs, their approach requires forward simulation of the system trajectories and/or a priori knowledge of bounds on the CBF dynamics. Alternatively, in \cite{Xiao2022Adaptive} an adaptation-based approach is introduced for guaranteed feasibility of a HO-CBF-QP control law provided that the parameters adapt sufficiently quickly. or bounded control authority \cite{Breeden2021InputConstraints, Xiao2022Adaptive}, these results do not generally extend to control systems seeking to satisfy multiple candidate CBF constraints. 

The problem of safe control design under multiple constraints is especially relevant in practical applications involving autonomous vehicles and mobile robots, where there may be liveness-based and/or spatiotemporal specifications in addition to safety constraints. In many cases, the joint satisfaction of safety and liveness constraints has been treated by synthesizing CBF-CLF-QP controllers \cite{Garg2022Fixed,Garg2019Control,Black2020Quadratic}, wherein CBF conditions are hard constraints and CLF conditions are soft constraints in the QP.
% , i.e., the CLF constraints are permitted to relax via some slack variable. 
As a class of Lyapunov-like functions, however, CBFs have also been used to enforce the satisfaction of tracking \cite{Shao2021Tracking} and spatiotemporal constraints using logic encodings like signal temporal logic (STL) \cite{Lindemann2019CBFSTL,Yang2020STL} and linear temporal logic (LTL) \cite{Srinivasan2021Control}. Whether there exists a control input capable of satisfying the full collection of CBF constraints, however, is very much still an open problem. 

Recent approaches to control design in the presence of multiple constraints have mainly circumvented the underlying problem by considering only one such constraint at a given time instance, either by assumption \cite{Cortez2022RobustMultiple} or construction in a non-smooth manner \cite{Glotfelter2017Nonsmooth,huang2020switched}, all of which may result in performance degradation (including undesirable oscillatory behavior). In contrast, the authors of \cite{Lindemann2019CBFSTL} and \cite{Machida2021ConsensusCBF} each propose smoothly synthesizing one candidate CBF for the joint satisfaction of multiple constraints, but notably make no attempt to validate their candidate function. Additional proposed solutions include learning suitable policies offline \cite{cai2021safe,Ma2021Generalized} and computing an a priori viability domain for multiple CBFs \cite{breeden2022compositions}, though as noted previously online methods provide better robustness to dynamic or unknown environments. Under such circumstances, in fact, various works have utilized parameter adaptation in an attempt to either learn \cite{Black2022Fixed,Lopez2021RaCBF} or compensate for \cite{Taylor2020Adaptive} unknown parameters in the system dynamics. 
In what follows, we explain how this paper uses parameter adaptation to progress toward a viable solution to control design under multiple spatiotemporal and input constraints.

\subsection{Summary of Contributions}
\begin{itemize}
    \item We formalize the notion of a consolidated CBF (C-CBF) by proposing a generic class of candidate functions that smoothly synthesizes any arbitrary number of time-varying constraint functions into one.
    \item We utilize a predictor-corrector interior point method to propose an optimization-based adaptation law for the weights of the C-CBF candidate's constituent functions. We show that the resulting weights render the C-CBF forward invariance condition viable at all times despite limited control authority, which serves to certify the C-CBF candidate as valid and permits the use of existing approaches to CBF-based control design for a class of nonlinear, control-affine systems.
    \item We illustrate the performance of our approach on a one-dimensional, academic example which requires the controller to recognize and avoid regions of the state space in which the specified constraints are satisfied, but where no control input is able to prevent the constraints from being violated in the future.
    \item Finally, we provide a comparative example utilizing a (car-like) robot with bicycle dynamics to highlight how our proposed controller solves a reach-avoid problem in which two baseline controllers fail, and highlight how the weight adaptation ensures that constraints are satisfied independent of their relative-degree.
\end{itemize}

% It is with this problem in mind that we propose an adaptive, consolidated CBF (C-CBF) based approach to constrained control design for classes of both single- and multi-agent systems. Constructed by smoothly synthesizing any arbitrary number of constraint functions into one, our C-CBF defines a new super-level set that under-approximates the intersection of its constituent sets arbitrarily closely (see Figure \ref{fig.adaptive_safe_set}). We further propose a parameter adaptation law for the weighting of the constituent functions, and prove that its use renders our C-CBF valid and the super-level set controlled invariant for the class of nonlinear, control-affine, multi-agent systems under consideration. And while various works have utilized parameter adaptation in the context of control for safety-critical systems, usually in an attempt to either learn \cite{Black2022Fixed,Lopez2021RaCBF} or compensate for \cite{Taylor2020Adaptive} unknown parameters in the system dynamics, our proposed adaptation law is the first to our knowledge to be used for the simultaneous verified satisfaction of multiple CBF constraints. To show the effectiveness of our proposed control formulation, we compare our approach to two state-of-the-art CBF-based controllers on an aerial robot reach-avoid problem and demonstrate how it handles constraints of arbitrary relative-degree. We then highlight its generalizability by studying a 25-agent decentralized control problem that requires each agent to satisfy 27 constraints at all times.

In our prior work \cite{black2022adaptation}, we introduced an adaptation law for rendering a similar form of consolidated CBF viable with respect to \textit{time-invariant} constraint functions under \textit{unbounded} control authority. Previously, the parameter adaptation law was taken to be the solution to a QP, as opposed to the predictor-corrector interior point method introduced in this paper. In addition, we extend our previous framework to handle spatiotemporal constraints (i.e., time-varying) and a class of input constraints. In this sense, the work herein has not been proposed to any other venue.

% Notably, the adaptation law in this paper is novel in that it is derived using interior point methods versus our previous work taking t 

% {\color{blue}
% It is worth noting that, in our prior work \cite{black2022adaptation}, we introduced a QP-based adaptation law for achieving such an objective for time-invariant constraint functions under \textit{unbounded} control authority. Importantly, however, our previous adaptation law was constructed without the use of interior point methods, the foundation on which the law proposed in this paper relies.
% }

The paper is organized as follows. Section \ref{sec.prelims} introduces some preliminaries, including time-varying convex optimization and constrained control design with CBFs. In Section \ref{sec.consolidated_cbf}, we introduce the form of our C-CBF and propose a parameter adaptation law for rendering it valid. Section \ref{sec.ccbf_control} introduces our proposed QP-based control law, the performance of which we demonstrate with two simulation examples in Section \ref{sec.case_studies}. In Section \ref{sec.conclusion} we conclude with final remarks and directions for future work.

% \subsection{Motivation}

% \subsection{Literature Review}

% It is worth further highlighting that while CBFs have most frequently been used as a tool for \texit{safe} control design \cite{ames2017control,etc}, 
% % in the sense that the constraint function $h_i$ implicitly defines a set of safe states $S_i$, 
% the fundamental machinery underpinning CBF-based control design is set invariance. Thus, a CBF may be used to enforce the satisfaction of any generic constraint, i.e., temporal specifications \cite{Yang2020Stl}, tracking requirements \cite{Shao2021Tracking}, etc. 

% \subsubsection{Convex Optimization}
% {\color{red} TO DO: Read this paper \cite{Rahili2017Convex}}

% In this paper \cite{Fazlyab2016Optimization}, the authors propose an approach to solve a constrained convex optimization problem using log-barrier methods -- this will be important in the development of our adaptation law.

%%%%%%%%%%%%%%%%%%%%%%%%%%%%%%%%%%%%%%%%%%%%%%%%%%%%%%%%%%
%************** Mathematical Preliminaries **************%
%%%%%%%%%%%%%%%%%%%%%%%%%%%%%%%%%%%%%%%%%%%%%%%%%%%%%%%%%%

\section{Mathematical Preliminaries}\label{sec.prelims}
We use the following notation throughout the paper. $\R$ denotes the set of real numbers and $\R_+$ the set of non-negative real numbers. The set of integers $\{1,\hdots,n\}$ is denoted $[n]$. Scalar quantities are denoted by lowercase variables ($t \in \R$), vectors are bold lowercase ($\bb{x} \in \R^n$), and matrices are bold uppercase ($\bb{M} \in \R^{n\times m}$). 
For a positive-definite matrix $\bb{M} \in \R^{n \times n}$ and scalar $a>0$, we write $\bb{M} \succeq a\bb{I}$ if all eigenvalues of $\bb{M}$ are greater than or equal to $a$, where $\bb{I}$ denotes identity of appropriate dimension. 
$\|\cdot\|$ represents the $L^2$ norm. A function $\alpha: \R \rightarrow \R$ belongs to extended class-$\mathcal{K}$ (expressed by $\alpha \in \mathcal{K}_\infty$) if $\alpha(0)=0$ and $\alpha$ is increasing on the interval $(-\infty,\infty)$. A function $\phi: \R \times \R \rightarrow \R$ belongs to class-$\mathcal{LL}$ (i.e., $\phi \in \mathcal{LL}$) if for each fixed $r$ (resp. $s$), the function $\phi(r,s)$ is decreasing with respect to $s$ (resp. $r$) and is such that $\phi (r,s) \rightarrow 0$ for $s \rightarrow \infty$ (resp. $r \rightarrow \infty$). The interior and boundary of a closed set $\mathcal{S}$ are denoted $\mathrm{Int}(\mathcal{S})$ and $\partial \mathcal{S}$ respectively. Given a multivariate function $V: \R \times \R^n$ denoted $V(t, \bb{x})$, let $\nabla _tV(t, \bb{x})$ be the derivative of $V$ with respect to $t \in \R$, $\nabla_{tt}V(t, \bb{x})$ its second derivative, $\nabla_{\bb{x}}V(t, \bb{x})$ the row vector of derivatives (i.e., gradient) of $V$ with respect to $\bb{x} \in \R^n$, and $\nabla_{\bb{xx}}V(t, \bb{x})$ its Hessian matrix.

% The Lie derivative of a function $V:\mathbb R^n\rightarrow \mathbb R$ along a vector field $f:\mathbb R^n\rightarrow\mathbb R^n$ at a point $x\in \mathbb R^n$ is denoted $L_fV(x) \triangleq \frac{\partial V}{\partial x} f(x)$. 

\subsection{Predictor-Corrector Interior Point Method}
Consider a time-varying, possibly non-convex optimization problem, a local solution trajectory to which takes the form
\begin{equation}\label{eq.opt_solution_trajectory}
\begin{aligned}\small
    \bb{y}^*(t) = \argmin_{\bb{y} \in \R^p} \; J(t, \bb{y}), \quad \textrm{s.t.} \;\; c_i(t, \bb{y}) \leq 0, \; \forall i \in [q],
\end{aligned}
\end{equation}
where $J: \R \times \R^p \rightarrow \R$ is a (uniformly strongly convex) objective function, the collection of $q \geq 0$ (not necessarily convex) inequality constraints is $c_i: \R \times \R^p \rightarrow \R$ for $i \in [q]$, and where all of $J$ and $c_i$ are twice continuously differentiable in $\bb{y}$ and piecewise continuously differentiable in $t$. Denote the feasible region of \eqref{eq.opt_solution_trajectory} as $\mathcal{Y}(t) = \{\bb{y} \in \R^p: c_i(t, \bb{y}) \leq 0, \forall i \in [q]\}$, and further assume that the feasible region has an interior point at all times, i.e., that $\mathrm{Int}(\mathcal{Y}(t)) \neq \emptyset$, $\forall t \geq 0$. Though $J$ is assumed to be uniformly strongly convex, without convexity of the inequality constraints it is not guaranteed that \eqref{eq.opt_solution_trajectory} is the global solution trajectory; nevertheless, barrier function methods may be employed to approximate a local solution \eqref{eq.opt_solution_trajectory} (see e.g., \cite[Sec. 9.4]{bazaraa2013nonlinear}). Consider such an approximate solution $\hat{\bb{y}}^*(t)$, constructed using Frisch's logarithmic barrier functions:
\begin{equation}\label{eq.approximate_solution}
\begin{aligned}\small
    \hat{\bb{y}}^*(t) = \argmin_{\bb{y} \in \R^p} \; \Psi(t, \bb{y}) \triangleq J(t, \bb{y}) - \frac{1}{s}\sum_{i=1}^q\log(-c_i(t, \bb{y})),
\end{aligned}
\end{equation}
where $s > 0$ is a barrier design parameter. In the limit, as $s \rightarrow \infty$, the function $-\frac{1}{s}\log(-c_i(t, \bb{y})) \rightarrow \mathbb{I}_{-}(c_i(t, \bb{y}))$, where $\mathbb{I}_{-}:\R \rightarrow \{0,\infty\}$ is defined such that $\mathbb{I}_{-}(u) = 0$ for $u \leq 0$ and $\mathbb{I}_{-}(u) = \infty$ for $u > 0$. In addition to facilitating this non-smooth function approximation, the barrier parameter may also convexify the augmented objective function $\Psi$. Though $\log(-c_i)$ may be non-convex, we assume (similar to e.g., \cite{bertsekas1979convexification}) that for sufficiently large $s$ the function $\Psi$ is strongly, uniformly convex, stated formally as follows.
\begin{Assumption}\label{ass.psi_convex}
    For all $\bb{y} \in \mathrm{Int}(\mathcal{Y}(t))$, $\forall t \geq 0$, there exist $s, a > 0$ such that
    \begin{equation*}
    \begin{aligned}\small
        \nabla_{\bb{yy}}\Psi(t, \bb{y}) \coloneqq \nabla_{\bb{yy}}J(t, \bb{y}) - \frac{1}{s}\sum_{i=1}^q\nabla_{\bb{yy}}\log(-c_i(t, \bb{y})) \succeq a \bb{I}.
        % \bb{\gamma}^T\left(\nabla_{\bb{yy}}J(\bb{z},t) - \frac{1}{s}\sum_{i=1}^q\nabla_{\bb{yy}}\log(c_i(\bb{y},t))\right)\bb{\gamma} \succeq a I
    \end{aligned}
    \end{equation*}
\end{Assumption}

In practice, for strongly uniformly convex $J$ the above may be satisfied by using sufficiently large $s$, which may be determined either online or iteratively offline. It follows then that $\nabla_{\bb{yy}}\Psi$ exists and is bounded, which allows us to review a version of \cite[Lemma 2]{Fazlyab2018Prediction} relating the approximate (locally) optimal trajectory $\hat{\bb{y}}^*(t)$ given by \eqref{eq.approximate_solution} to the solution of a continuous-time dynamical system.
\begin{Lemma}\label{lem.feasible_optimization_solution}
    Let $\hat{\bb{y}}^*(t)$ be defined by \eqref{eq.approximate_solution} and $\bb{z}(t)$ be the solution to the following ordinary differential equation,
    \begin{equation}\label{eq.zdot_ode_optimization}
    \begin{aligned}
        \dot{\bb{z}} &= -\nabla_{\bb{zz}}^{-1}\Psi(t, \bb{z})\big[\bb{P}\nabla_{\bb{z}}\Psi(t, \bb{z}) + \nabla_{\bb{z}t}\Psi(\bb{z},t)\big],\\
        \bb{z}(0) &= \bb{z}_0 \in \mathcal{Y}(0),
    \end{aligned}
    \end{equation}
    where $\bb{P} \in \R^{p \times p}$ is a positive-definite gain matrix satisfying $\bb{P} \succeq b \bb{I}$ for some $b > 0$, and $\Psi: \R \times \R^{p}$ is defined as in \eqref{eq.approximate_solution} and known to satisfy Assumption \ref{ass.psi_convex}. Then, $\bb{z}(t) \in \mathcal{Y}(t)$ for all $t \geq 0$ and
    \begin{equation}\nonumber
        \|\bb{z}(t) - \hat{\bb{y}}^*(t)\| \leq Ce^{-b t},
    \end{equation}
    where $0 \leq C \coloneqq \frac{1}{a}\|\nabla_{\bb{z}}\Psi(\bb{z}(0),0)\| < \infty$.
\end{Lemma}

The above result implies that the solution $\bb{z}(t)$ to \eqref{eq.zdot_ode_optimization} approaches the approximate optimal trajectory $\hat{\bb{y}}^*(t)$ asymptotically while remaining in the feasible region $\mathcal{Y}(t)$ at all times. Intuitively, the dynamics of \eqref{eq.zdot_ode_optimization} consist of a prediction term
\begin{equation*}
    \dot{\bb{z}}_p = -\nabla_{\bb{zz}}^{-1}\Psi(t, \bb{z})\nabla_{\bb{z}t}\Psi(t, \bb{z}),
\end{equation*}
and a Newton-like correction term
\begin{equation*}
    \dot{\bb{z}}_c = -\nabla_{\bb{zz}}^{-1}\Psi(t, \bb{z})\bb{P}\nabla_{\bb{z}}\Psi(t, \bb{z}),
\end{equation*}
which together ensure that $\bb{z}(t) \in \mathrm{Int}(\mathcal{Y}(t))$ at all times (hence predictor-corrector interior point method).
% such that the solution to \eqref{eq.zdot_ode_optimization} asymptotically tracks the approximate optimal solution trajectory $\hat{y}^*(t)$ while adhering to the feasible region $\mathcal{Y}(t)$ at all times.. 
We require the findings from Lemma \ref{lem.feasible_optimization_solution} in the proof of one of our main results on parameter adaptation for constrained control design, the premises of which we will now review.

\subsection{Constrained Control Design}\label{subsec.constrained_control}

In the remainder of this paper we consider a class of nonlinear, control-affine dynamical systems of the form
\begin{equation}\label{eq.single_agent_system}
    \dot{\bb{x}} = f(\bb{x}(t)) + g(\bb{x}(t))\bb{u}(t), \quad \bb{x}(0) = \bb{x}_{0}
\end{equation}
where $t \in \mathcal{T} = [t_0, \infty)$ represents time, $\bb{x} \in \R^n$ and $\bb{u} \in \mathcal{U} \subset \R^m$ are the state and control input vectors, with $\mathcal{U}$ the input constraint set, and where $f: \R^n \rightarrow \R^n$ and $g: \R^n \rightarrow \R^{n \times m}$ are known and locally Lipschitz. The input constraint set is
\begin{equation}\label{eq.input_constraint_set}
    \mathcal{U} = \{\bb{u} \in \R^m \mid -\bb{u}_{max} \leq \bb{u} \leq \bb{u}_{max}\},
\end{equation}
where the inequalities are interpreted element-wise and $\bb{u}_{max} = (\bar{u}_{1},\cdots,\bar{u}_{m})$ with $0 < \bar{u}_j < \infty$ for $j \in [m]$. We assume the availability of a control input $\bb{u}: \mathcal{T} \rightarrow \mathcal{U}$ that yields a unique solution to \eqref{eq.single_agent_system} for all $t \in \mathcal{T}$.

Consider one of $c \geq 1$ (possibly time-varying) twice continuously differentiable constraint functions $h_i: \mathcal{T} \times \R^n \rightarrow \R$,
% in $\bb{x}$ and piecewise continuously differentiable in $t$
which may encode a safety, performance, and/or specification based constraint, and let the set of states over which this particular constraint is satisfied at a given time $t$ be
\begin{equation}\label{eq.safe_set}
    \mathcal{S}_i(t) = \{\bb{x} \in \R^n \mid h_i(t, \bb{x}) \geq 0\}.
\end{equation}
We denote the set of all constraint functions as 
\begin{equation}\label{eq.constraint_function_set}
    \mathcal{H} = \{h_1,\cdots,h_c\},
\end{equation}
and in this paper focus on confining the state to the intersection of all individual constraint sets, hereafter referred to as the complete constraint set:
\begin{equation}\label{eq.complete_constraint_set}
    \mathcal{S}(t) = \bigcap_{i=1}^{c}\mathcal{S}_i(t).
\end{equation}
We assume that the complete constraint set possesses an interior point, i.e $\mathrm{Int}(\mathcal{S}(t)) \neq \emptyset$, $\forall t \in \mathcal{T}$, and thus seek to design a controller that renders a subset $\mathcal{D}(t) \subseteq \mathcal{S}(t)$ \textit{forward-invariant} with respect to the system \eqref{eq.single_agent_system}, i.e., to ensure that $\bb{x}(0) \in \mathcal{D}(0) \subseteq \mathcal{S}(0) \implies \bb{x}(t) \in \mathcal{D}(t) \subseteq \mathcal{S}(t)$, $\forall t \in \mathcal{T}$. 
% For many safety-critical applications it is reasonable to assume that the constraint functions are compatible, e.g., a vehicle is not required be on both sides of a wall at once. Under a mix of safety-critical and performance-driven constraints, however, it may be necessary to relax a performance constraint function, e.g., $h^r_i(t, \bb{x}) = h_i(t, \bb{x}) + \epsilon$ for some $\epsilon > 0$, in order to preserve the interior of $\mathcal{S}$. We address this issue in Section \ref{subsec.relaxation}.
Of note is that we require that the functions $h_i$ be twice continuously differentiable, a stricter (though often reasonable) condition than the more typical continuous differentiability requirement, for reasons pertaining to our parameter adaptation law in Section \ref{sec.consolidated_cbf}.

For the remainder of this section we will introduce preliminaries specific to the case of one constraint function ($c=1$). Under such circumstances, one approach to controlled set invariance is to use CBFs in the control design. We now review the notion of a time-varying CBF inspired by \cite{breeden2021robust,Lindemann2019CBFSTL}, and highlight how it is used to enforce forward-invariance.
\begin{Definition}\label{def.cbf}
    Given a time-varying set $\mathcal{S}_i(t)$, $t \in \mathcal{T}$, defined by \eqref{eq.safe_set} for a twice continuously differentiable function $h_i: \mathcal{T} \times \R^n \rightarrow \R$, the function $h_i$ is a \textbf{control barrier function} (CBF) for the system \eqref{eq.single_agent_system} with respect to $\mathcal{S}_i(t)$ if there exists a locally Lipschitz continuous function $\alpha \in \mathcal{K}_\infty$ such that, for all $\bb{x} \in \mathcal{S}_i(t)$, $t \in \mathcal{T}$,
    \begin{equation}\label{eq.tv_cbf_condition}
        \sup_{\bb{u} \in \mathcal{U}}\;\dot{h}_i(t, \bb{x}, \bb{u}) \geq -\alpha(h_i(t, \bb{x})).
    \end{equation}
\end{Definition}
We refer to \eqref{eq.tv_cbf_condition} as the CBF condition, and note that if the function $h_i$ has relative-degree\footnote{A function $p: \R_+ \times \R^n \rightarrow \R$ is said to be of relative-degree $r$ with respect to the dynamics \eqref{eq.single_agent_system} if $r$ is the number of times $p$ must be differentiated before one of the control inputs $u$ appear explicitly.} $r > 1$ with respect to the system \eqref{eq.single_agent_system} then the control input $\bb{u}$ has no effect on \eqref{eq.tv_cbf_condition}. In many cases this deficiency can be resolved by deriving high-order CBFs (see e.g., \cite{Breeden2021InputConstraints,Xiao2019HOCBF} for details), though this may require careful parameter design \cite{nguyen2016exponential}. The method we propose in Section \ref{sec.consolidated_cbf}, however, requires no assumptions about the relative-degree of the constraint functions under consideration, a fact that we highlight in our second numerical case study in Section \ref{sec.case_studies}. For a particular constraint function $h_i$, we refer to the set of control actions satisfying its CBF condition \eqref{eq.tv_cbf_condition} as its CBF control set and denote it by
\begin{align}
    \mathcal{U}_{h_i}(t, \alpha) = \big\{&\bb{u} \in \mathcal{U}: \label{eq.admissible_control_set} \\
    &\dot{h}_i(t, \bb{x}(t), \bb{u}) +\alpha\big(h_i(t, \bb{x}(t))\big) \geq 0\big\}, \nonumber
\end{align}
noting that $\mathcal{U}_{h_i}$ is parameterized by the set function $\alpha \in \mathcal{K}_{\infty}$. In the following result inspired by \cite{Lindemann2019CBFSTL}, it is shown that controls belonging to the CBF control set $\mathcal{U}_{h_i}$ render the constraint set $\mathcal{S}_i(t)$ forward-invariant.

\begin{Theorem}\label{thm.cbf_invariance}
    If $h_i$ is a control barrier function for the system \eqref{eq.single_agent_system} with respect to the set $\mathcal{S}_i(t)$ defined by \eqref{eq.safe_set}, then $\bb{x}(0) \in \mathcal{S}_i(0) \implies \bb{x}(t) \in \mathcal{S}_i(t)$ for all $t \in \mathcal{T}$, i.e., $\mathcal{S}_i(t)$ is forward-invariant.
\end{Theorem}
\begin{proof}
    Consider an absolutely continuous function $\eta: \mathcal{T} \rightarrow \R$ and locally Lipschitz function $\gamma \in \mathcal{K}_\infty$. By \cite[Lem. 2]{Glotfelter2017Nonsmooth} it is true that if $\dot\eta(t) \geq -\gamma(\eta(t))$ for every $t \in \mathcal{T}$ and $\eta(0) \geq 0$, then $\eta(t) \geq 0$ for all $t \in \mathcal{T}$. With assumed unique solutions to \eqref{eq.single_agent_system} and taking $\eta(t) = h_i(t, \bb{x}(t))$ and $\gamma = \alpha$, it follows that $h_i(t, \bb{x}(t)) \geq 0$ for all $t \in \mathcal{T}$ and thus $\mathcal{S}_i(t)$ is forward-invariant.
\end{proof}

Note that $h_i$ being a CBF implies that there exists $\alpha \in \mathcal{K}_\infty$ such that the CBF control set is non-empty for all time, i.e., $\mathcal{U}_{h_i}(t, \alpha) \neq \emptyset$, $\forall t \in \mathcal{T}$. Under such circumstances, for the class of control-affine systems described by \eqref{eq.single_agent_system} the CBF condition \eqref{eq.tv_cbf_condition} may be included as a linear constraint in the following quadratic program (QP) based control law (see e.g., \cite{Black2022ffcbf,ames2017control}):
\begin{subequations}\label{eq.cbf_qp_controller}
\begin{align}
    \bb{u}^*(t) = \argmin_{\bb{u} \in \mathcal{U}} \; &\frac{1}{2}\|\bb{u}-\bb{u}_0(t, \bb{x})\|^2 \label{subeq.cbf_qp_objective}\\
    \textrm{s.t.} \nonumber \\
    \frac{\partial h_i}{\partial t} + \frac{\partial h_i}{\partial \bb{x}}f(\bb{x}) &+ \frac{\partial h_i}{\partial \bb{x}}g(\bb{x})\bb{u} \geq -\alpha\big(h_i(t, \bb{x})\big), \label{subeq.cbf_qp_constraints}
\end{align}
\end{subequations}
where \eqref{subeq.cbf_qp_objective} seeks to produce a solution $\bb{u}^* \in \mathcal{U}$ that deviates minimally from some desired input $\bb{u}_0: \mathcal{T} \times \R^n \rightarrow \R^m$, and \eqref{subeq.cbf_qp_constraints} encodes the CBF condition \eqref{eq.tv_cbf_condition}. For systems without input constraints (i.e., $\mathcal{U} = \R^m$), other works (e.g., \cite{Jankovic2022Multi}) highlight that the constraint function $h_i$ is a CBF if there exists a function $\alpha \in \mathcal{K}_\infty$ satisfying
\begin{equation}\label{eq.alternate_tv_cbf_condition_unbounded}
\begin{aligned}
    \frac{\partial h_i}{\partial \bb{x}}g(\bb{x}) = \mathbf{0} &\implies \\
    &\frac{\partial h_i}{\partial t} + \frac{\partial h_i}{\partial \bb{x}}f(\bb{x}) + \alpha(h_i(t, \bb{x})) > 0
\end{aligned}
\end{equation}
for all $\bb{x}(t) \in \mathcal{S}_i(t)$, $t \in \mathcal{T}$. For the class of systems under consideration in this paper (with $\mathcal{U} \neq \R^m$), the analogous requirement is that 
\begin{equation}\label{eq.alternate_tv_cbf_condition_bounded}
\begin{aligned}
    \sup_{\bb{u} \in \mathcal{U}}\frac{\partial h_i}{\partial \bb{x}}g(\bb{x})\bb{u} =& \delta \implies \\
    &\frac{\partial h_i}{\partial t} + \frac{\partial h_i}{\partial \bb{x}}f(\bb{x}) + \alpha(h_i(t, \bb{x})) > -\delta
\end{aligned}
\end{equation}
for all $\bb{x}(t) \in \mathcal{S}_i(t)$, $t \in \mathcal{T}$, $\delta \geq 0$. Here, $\delta$ serves to define the strongest drift to the CBF dynamics which the controller can overcome to satisfy the CBF condition \eqref{eq.tv_cbf_condition}. It is worth highlighting that both cases \eqref{eq.alternate_tv_cbf_condition_unbounded} and \eqref{eq.alternate_tv_cbf_condition_bounded} are equivalent to the condition that the CBF control set $\mathcal{U}_{h_i}$ is never empty. While for systems with unlimited control authority it may be straightforward to show that $\mathcal{U}_{h_i} \neq \emptyset$ (see e.g., the collision avoidance constraint between two vehicles in \cite{Jankovic2021Collision}), verifying that \eqref{eq.alternate_tv_cbf_condition_bounded} is true a priori for a given constraint function is generally much more difficult. Unless \eqref{eq.alternate_tv_cbf_condition_bounded} is proven to hold, the constraint function $h_i$ is only a \textit{candidate} CBF for \eqref{eq.single_agent_system} with respect to $\mathcal{S}_i(t)$, and thus may or may not render $\mathcal{S}_i(t)$ forward-invariant. Recent works have addressed the problem of synthesizing a \textit{valid} CBF via offline analysis \cite{clark2021verification}, 
% \cite{srinivasan2020synthesis}, 
by construction \cite{Breeden2021InputConstraints}, inclusion of additional constraints in the QP controller \cite{xiao2022sufficient}, or through online adaptation \cite{Xiao2022Adaptive}. The scopes of these solutions, however, are limited to only one time-invariant constraint function ($c = 1$) and therefore do not extend to situations in which the state must remain within a complete constraint set of the form \eqref{eq.complete_constraint_set} comprised of multiple, time-varying constituent constraint functions.

\subsection{Problem Statements}

% In this paper, we seek to contain the state of system \eqref{eq.single_agent_system} within the complete constraint set defined by \eqref{eq.complete_constraint_set} for a collection of $c \geq 1$ constraint functions. This challenge is more representative of real-world control problems where many safety and performance specifications may be entwined. 
% And while each constraint may be encoded as a candidate CBF, 
The main open problem in constrained control design under multiple state and input constraints is in certifying that all constraints will satisfied \textit{jointly} in perpetuity. We bridge this gap by 1) synthesizing one \textit{consolidated} CBF (C-CBF) candidate from the collection of constraints and then 2) certifying the C-CBF candidate as \textit{valid} via online parameter adaptation. 
% With this approach we may then use a controller of the form \eqref{eq.cbf_qp_controller} directly, which results in certifiable constrained control design. 
% We now formally state the first problem under consideration.
\begin{Problem}\label{prob.consolidated_cbf}
    Consider a collection of $c \geq 1$ twice continuously differentiable constraint functions $h_i: \mathcal{T} \times \R^n \rightarrow \R$ corresponding to constraint sets $\mathcal{S}_i(t)$ given by \eqref{eq.safe_set}. Design a new constraint function $H: \mathcal{T} \times \R_+^c \times \R^n \rightarrow \R$ with constituent constraint function weights $\bb{w} = (w_1,\cdots,w_c) \in \R_+^c$ such that the set 
    \begin{equation}\label{eq.ccbf_zero_superlevel_set}
        \mathcal{D}(t, \bb{w}) = \left\{\bb{x} \in \R^n \mid H(t, \bb{w}, \bb{x}) \geq 0\right\}
    \end{equation}
     satisfies $\mathcal{D}(t, \bb{w}) \subseteq \mathcal{S}(t)$ for the complete constraint set given by \eqref{eq.complete_constraint_set}, $\forall t \in \mathcal{T}$.
\end{Problem}

We refer to the new constraint function $H$ as a C-CBF candidate, and seek to guarantee that we can render its zero super-level set $\mathcal{D}(t, \bb{w})$ forward-invariant in the presence of input constraints of the form \eqref{eq.input_constraint_set}. Consolidating many separate constraint functions $h_i$ into one function $H$ means that each weight $w_i$ affects the relative importance of $h_i$ in $H$ (and by extension on the set $\mathcal{D}$). It also allows us, however, to leverage existing CBF results for certifiable constrained control design. Notably, we may then use the class of QP based controllers \eqref{eq.cbf_qp_controller} with the new function $H$ to render $\mathcal D$ forward-invariant and thus to ensure that the state remains within the complete constraint set $\mathcal S$. Proving that $H$ is a valid C-CBF, however, may be difficult when using static constituent function weights $\bb{w}$, i.e., weights chosen a priori with $\dot{\bb{w}} \equiv \mathbf{0}_{c \times 1}$. We therefore seek to design a weight adaptation law $\dot{\bb{w}} = \omega(t, \bb{w}, \bb{x}, \bb{u})$, $\omega: \mathcal{T} \times \R_+^c \times \R^n \times \mathcal{U} \rightarrow \R^c$, that varies the weights $\bb{w}$ online and, in doing so, renders the C-CBF candidate $H$ valid.
With $\bb{w}$ then dependent on time $t$, in the remainder we omit the dependence of $\mathcal{D}$ on the weights $\bb{w}$ for conciseness.
To render the C-CBF candidate $H$ valid is to ensure that
\begin{equation}\label{eq.c-cbf-condition}
    \sup_{\bb{u} \in \mathcal{U}}\dot H(t, \bb{w}, \dot{\bb{w}}, \bb{x}, \bb{u}) \geq -\alpha\big(H(t, \bb{w}, \bb{x})\big),
\end{equation}
holds for all $\bb{x}(t) \in \mathcal{D}(t)$, $\forall t \in \mathcal{T}$. Notice that the above is to the C-CBF candidate $H$ as \eqref{eq.tv_cbf_condition} is to a CBF candidate $h_i$. Thus, we refer to \eqref{eq.c-cbf-condition} as the C-CBF condition. While we have established that the CBF condition \eqref{eq.tv_cbf_condition} is affine in the control input $\bb{u}$ for the class of systems described by \eqref{eq.single_agent_system}, the C-CBF condition \eqref{eq.c-cbf-condition} is only affine in $\bb{u}$ if the adaptation law $\dot{\bb{w}}$ is affine in (or independent of) $\bb{u}$, i.e., if $\dot{\bb{w}}$ takes the form
\begin{equation}\label{eq.control_affine_adaptation}
    \dot{\bb{w}} = \omega(t, \bb{w}, \bb{x}, \bb{u}) \coloneqq \mu(t, \bb{w}, \bb{x}) + \nu(t, \bb{w}, \bb{x})\bb{u},
\end{equation}
for some functions $\mu: \mathcal{T} \times \R_+^c \times \R^n \rightarrow \R^c$ and $\nu: \mathcal{T} \times \R_+^c \times \R^n \rightarrow \R^{c \times m}$. With $\dot{\bb{w}}$ taking this form, we may directly include the C-CBF condition \eqref{eq.c-cbf-condition} as a linear constraint in a QP-based control law of the form \eqref{eq.cbf_qp_controller}. As such, we now formally introduce the second problem under consideration.
\begin{Problem}\label{prob.adaptation}
    Given a function $\alpha \in \mathcal{K}_\infty$, a consolidated control barrier function candidate $H$, and its corresponding CBF control set 
    \begin{align}
        \mathcal{U}_{H}(t, \alpha) = \big\{&\bb{u} \in \mathcal{U}: \label{eq.ccbf_control_set} \\
        &\dot{H}(t, \bb{w}, \dot{\bb{w}}, \bb{x}, \bb{u}) +\alpha\big(H(t, \bb{w}, \bb{x})\big) \geq 0\big\}, \nonumber
    \end{align}
    design a control-affine weight adaptation law of the form \eqref{eq.control_affine_adaptation} such that $\mathcal{U}_H(t, \alpha) \neq \emptyset$ for all $t \in \mathcal{T}$, i.e., such that
    \begin{align}
        \sup_{\bb{u} \in \mathcal{U}}\left(\frac{\partial H}{\partial \bb{x}}g(\bb{x}) + \frac{\partial H}{\partial \bb{w}}\nu(t, \bb{w}, \bb{x})\right)\bb{u} &= \delta \implies \label{eq.c-cbf-delta-condition} \\
        \frac{\partial H}{\partial t} + \frac{\partial H}{\partial \bb{x}}f(\bb{x}) + \frac{\partial H}{\partial \bb{w}}\mu(t, \bb{w}, \bb{x}) &+ \alpha\big(H(t, \bb{w}, \bb{x})\big) > -\delta, \nonumber
    \end{align}
    % \begin{align}
    %     \sup_{\bb{u} \in \mathcal{U}}\frac{\partial H}{\partial \bb{x}}g(\bb{x})\bb{u} &= \delta \implies \label{eq.c-cbf-delta-condition} \\
    %     \frac{\partial H}{\partial t} &+ \frac{\partial H}{\partial \bb{x}}f(\bb{x}) + \frac{\partial H}{\partial \bb{w}}\dot{\bb{w}} + \alpha\big(H(t, \bb{w}, \bb{x})\big) > -\delta, \nonumber
    % \end{align}
    for all $\bb{x}(t) \in \mathcal{D}(t)$, $\forall t \in \mathcal{T}$, $\delta \geq 0$.
\end{Problem}

Notice that \eqref{eq.c-cbf-delta-condition} constitutes a sufficient condition for the satisfaction of the C-CBF condition \eqref{eq.c-cbf-condition}, and therefore ensures that the CBF control set $\mathcal{U}_H(t,\alpha)$ is never empty. 
% In the ensuing section, we introduce our solutions to Problem \ref{prob.consolidated_cbf}, a generic form for the consolidated CBF candidate, and Problem \ref{prob.adaptation}, an optimization-based parameter adaptation law.

%%%%%%%%%%%%%%%%%%%%%%%%%%%%%%%%%%%%%%%%%%%%%%%%%%%%%%%%%%
%******************* Consolidated CBF *******************%
%%%%%%%%%%%%%%%%%%%%%%%%%%%%%%%%%%%%%%%%%%%%%%%%%%%%%%%%%%
\section{Adaptation for Consolidated CBF Synthesis}\label{sec.consolidated_cbf}
In this section, we first introduce our proposed solution to Problem \ref{prob.consolidated_cbf}, a consolidated control barrier function (C-CBF) candidate that smoothly synthesizes multiple constraint functions into one whose zero super-level set under-approximates the complete constraint set, and then we design a parameter adaptation law which renders the candidate C-CBF valid for constrained control design.

\subsection{Consolidated CBFs}
\begin{figure}[!ht]
    \centering
        \includegraphics[clip,width=1\linewidth]{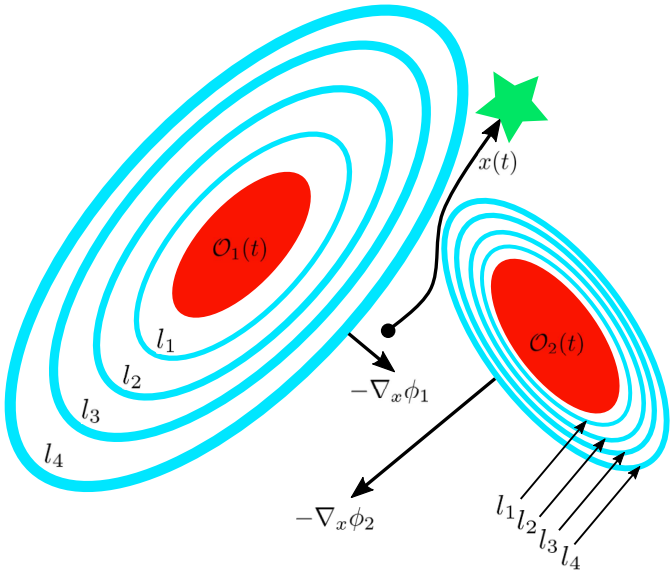}
    \caption{A contrived example of a system seeking to reach a target location (green star) while adapting weights $\bb{w}$ in the presence of two elliptical obstacles of identical size and shape whose constraint functions $h_1$ and $h_2$ define exclusion sets $\mathcal{O}_1(t) = \R^n \setminus \mathcal{S}_1(t)$ and $\mathcal{O}_2(t) = \R^n \setminus \mathcal{S}_2(t)$. Here, the system gives greater relative importance to obstacle one (i.e., $w_1 < w_2$) and therefore the function $H$ views obstacle one as having a shallow level curve topography (depicted as level curves $l_1<l_2<l_3<l_4$ spaced further apart with smaller gradient magnitude $\nabla_x\phi_1$), whereas it sees obstacle two as having steep gradients (level curves close together, with larger $\nabla_x\phi_2$) near the obstacle. This allows the state to more closely approach obstacle two before deciding to take evasive action.}\label{fig.level_set_example}
\end{figure}

% Denote by $\bb{h}(t, \bb{x})$ and $\bb{w}$ respectively the vectors of $c\geq 1$ constituent constraint functions evaluated at the state and their corresponding weights, i.e.,
% \begin{align}
%     \bb{h}(t, \bb{x}) &= \big(h_1(t, \bb{x}), \cdots, h_c(t, \bb{x})\big) \in \R^c, \nonumber \\
%     \bb{w} &= \big(w_1, \cdots, w_c\big) \in \R_{+}^c. \nonumber
% \end{align}
The form of our C-CBF candidate $H: \mathcal{T} \times \R^c_+ \times \R^n \rightarrow \R$ is
\begin{equation}\label{eq.consolidated_cbf}
    H(t, \bb{w}, \bb{x}) = 1 - \sum_{i=1}^c\phi\Big(h_i(t, \bb{x}), w_i\Big),
\end{equation}
where $\phi: \R_+ \times \R_+ \rightarrow \R_+$ belongs to class-$\mathcal{L}\mathcal{L}$, is twice continuously differentiable, and satisfies $\phi(r,0)=\phi(0,s)=\phi(0,0)=1$. Note that the specified domain considers only positive values for $w_i$ and $h_i$, a reasonable choice since $\bb{x}(t) \notin \mathcal{D}(t)$ if some $h_i(t,\bb{x}(t)) \leq 0$ with $w_i \geq 0$ or if some $w_i \leq 0$ with $h_i(t,\bb{x}(t)) \geq 0$. Both the decaying exponential function $\phi(r,s)=e^{-rs}$ and the class of reciprocal functions of the form $\phi(r,s)=v / (rs + v)$, $v > 0$, for example, satisfy the above requirements over the admissible domain. The function $\phi$ thus encodes that the zero super-level set of $H$ is a subset of the complete constraint set $\mathcal{S}(t)$, i.e.,
\begin{equation}\nonumber
    \mathcal{D}(t) = \{\bb{x} \in \R^n \mid H(t, \bb{w}(t),\bb{x}) \geq 0\} \subset \mathcal{S}(t), \; \forall t \in \mathcal{T},
\end{equation}
provided that $\bb{w} \in \R^c_{+}$. Note that the weights behave as a shape parameter for the function $\phi$, and are neither used to construct a weighted average nor required to sum to any specific amount. Instead, observe that larger weights $\bb{w}$ allow $\mathcal{D}(t)$ to more closely approach $\mathcal{S}(t)$, as for a given state it follows that $\phi(h_i(t, \bb{x}), w_i^l) < \phi(h_i(t, \bb{x}), w_i^s)$ for $w_i^l > w_i^s > 0$. This means that a higher weight $w_i$ confers smaller relative importance of its associated constraint function $h_i$ to the C-CBF candidate $H$. For example, in Figure \ref{fig.level_set_example} the effect of two different weights on the level sets of two obstacle avoidance constraint functions is depicted. In this case, obstacle two has smaller relative importance (i.e., $w_2 > w_1$) and thus the system takes actions toward obstacle two until it reaches a region where the level set values decrease very quickly to zero, at which point it makes an evasive maneuver. Thus, the system effectively ignores obstacle two until it must take action to avoid it.

We have then that $H$ defined by \eqref{eq.consolidated_cbf} solves Problem \ref{prob.consolidated_cbf} and is referred to as a C-CBF \textit{candidate}. Given an arbitrary weight vector $\bb{w}$, however, it may not be true that $H$ is a \textit{valid} CBF for \eqref{eq.single_agent_system} with respect to the set $\mathcal{D}(t)$. Without adapting the weights as the system moves through the state space the C-CBF condition \eqref{eq.c-cbf-delta-condition} may be violated.
% without adapting $\bb{w}$, i.e., for $\dot{\bb{w}} \equiv \mathbf{0}_{c \times 1}$. For example, if $\dot{\bb{w}} = \mu(t, \bb{w}, \bb{x}) + \nu(t, \bb{w}, \bb{x})\bb{u} \equiv \mathbf{0}_{c \times 1}$ it may occur that \eqref{eq.c-cbf-delta-condition} is violated, i.e., that
% \begin{equation*}\small
%     \sup_{\bb{u} \in \mathcal{U}}\frac{\partial H}{\partial \bb{x}}g(\bb{x})\bb{u} &= -\inf_{\bb{u} \in \mathcal{U}}\sum_{i=1}^c\frac{\partial \phi(h_i(t, \bb{x}), w_i)}{\partial h_i}\frac{\partial h_i}{\partial \bb{x}}g(\bb{x})\bb{u} < \delta,
% \end{equation*}
% when
% \begin{equation*}\small
%     \frac{\partial H}{\partial t} + \frac{\partial H}{\partial \bb{x}}f(\bb{x}) + \alpha(H(t, \bb{w}, \bb{x})) \leq -\delta,
% \end{equation*}
% for some $\delta>0$, in which case the C-CBF condition \eqref{eq.c-cbf-condition} would fail to hold. 
In what follows we introduce an online parameter adaptation law to enforce \eqref{eq.c-cbf-delta-condition} and, in doing so, to solve Problem \ref{prob.adaptation}.

\subsection{Adaptation for Online Verification}

Consider that in addition to the C-CBF weights $\bb{w}$ being used to enforce the condition \eqref{eq.c-cbf-condition}, it may be desirable to choose weights that are optimal with respect to some objective function $J: \mathcal{T} \times \R^c \times \R^n \rightarrow \R$. As such, we propose using a predictor-corrector interior point method for our adaptation law. 

Consider the following optimal solution trajectory:
\begin{equation}\label{eq.w_solution_trajectory}
\begin{aligned}
    \bb{w}^*(t) &= \argmin_{\bb{w} \in \R^c} \; J(t, \bb{w}, \bb{x}), \;\; \textrm{s.t.} \;\; \bb{w} \in \mathcal{W}(t),
\end{aligned}
\end{equation}
where $J$ is strongly uniformly convex with respect to $\bb{w}$, twice continuously differentiable in both $\bb{w}$ and $\bb{x}$, and piecewise continuously differentiable in $t$. The feasible region is
\begin{equation}\label{eq.w_feasible_region}
    \mathcal{W}(t) = \left\{\bb{w} \in \R^c \mid b_j(t, \bb{w}, \bb{x}) \leq 0, \forall j \in [2c+1]\right\},
\end{equation}
where the first $c$ constraint functions $b_j: \mathcal{T} \times \R^c \times \R^n \rightarrow \R$ enforce that each $w_j$ remains above a user-specified minimum threshold $w_{min} > 0$ to ensure that $\bb{w}(t) \in \R^c_+, \forall t \in \mathcal{T}$, i.e.,
\begin{equation}\label{eq.b_kmin}
    b_j(t, \bb{w},\bb{x}) = w_{min} - w_j, \; \forall j \in [c],
\end{equation}
the second $c$ constraint functions $b_j: \mathcal{T} \times \R^c \times \R^n \rightarrow \R$ enforce that each $w_j$ remains below a user-specified maximum threshold $w_{max} < \infty$ to ensure that $\bb{w}(t)$ is well-defined, i.e.,
\begin{equation}\label{eq.b_kmax}
    b_{j+c}(t, \bb{w},\bb{x}) = w_j - w_{max}, \; \forall j \in [c],
\end{equation}
and the final constraint function $b_{2c+1}: \mathcal{T} \times \R^c \times \R^n \rightarrow \R$ enforces that the C-CBF condition \eqref{eq.c-cbf-condition} is satisfied, i.e., namely that
\begin{align}
    \beta(t, \bb{w}, \bb{x}) \coloneqq &-\frac{\partial H}{\partial t} - \frac{\partial H}{\partial \bb{x}}f(\bb{x}) - \frac{\partial H}{\partial \bb{w}}\dot{\bb{w}} \label{eq.b_ccbf_condition_unfiltered} \\
    &- \alpha(H(t, \bb{w}, \bb{x})) - \sup_{\bb{u} \in \mathcal{U}}\left[\frac{\partial H}{\partial \bb{x}}g(\bb{x})\bb{u}\right] \leq 0. \nonumber
\end{align}
For the moment, let us consider that $b_{2c+1} = \beta$.
Recognizing that the optimal solution trajectory \eqref{eq.w_solution_trajectory} is of the form \eqref{eq.opt_solution_trajectory}, we similarly consider the approximate optimal solution trajectory
\begin{equation}\label{eq.what_solution_trajectory}
    \hat{\bb{w}}^*(t) = \argmin_{\bb{w} \in \R^c} \; \Phi(t, \bb{w}, \bb{x}),
\end{equation}
where $\Phi$ is analogous to $\Psi$ in \eqref{eq.approximate_solution} and is defined by
\begin{equation}\label{eq.Phi_func}
    \Phi(t, \bb{w}, \bb{x}) \coloneqq J(t, \bb{w}, \bb{x}) - \frac{1}{s}\sum_{j=1}^{2c+1}\log\big(-b_j(t, \bb{w}, \bb{x})\big),
\end{equation}
for barrier parameter $s>0$. We assume that Assumption \ref{ass.psi_convex} holds for $\Phi$ i.e., that there exists $s>0$ rendering $\Phi$ strongly uniformly convex with respect to $\bb{w}$. As stated previously, this is reasonable in practice given a sufficiently large choice of $s$. We therefore propose to adapt the weights $\bb{w}$ according to a control-affine dynamical system of the form \eqref{eq.zdot_ode_optimization}:
\begin{equation}\label{eq.w_adaptation_law}
    \dot{\bb{w}}= \mu(t, \bb{w}, \bb{x}) + \nu(t, \bb{w}, \bb{x})\bb{u},
\end{equation}
where $\mu: \mathcal{T} \times \R^c \times \R^n \rightarrow \R^c$ and $\nu: \mathcal{T} \times \R^c \times \R^n \rightarrow \R^{c \times m}$ are defined by
\begin{align}
    \mu(t, \bb{w}, \bb{x}) &= -\nabla_{\bb{ww}}^{-1}\Phi\big(\bb{P}\nabla_{\bb{w}}\Phi + \nabla_{\bb{wx}}\Phi f(\bb{x}) + \nabla_{\bb{w}t}\Phi\big), \label{eq.ka_dot} \\
    \nu(t, \bb{w}, \bb{x}) &= -\nabla_{\bb{ww}}^{-1}\Phi\nabla_{\bb{wx}}\Phi g(\bb{x}), \label{eq.kb_dot}
\end{align}
for which $\bb{P} \in \R^{c \times c}$ is a positive-definite gain matrix. Note that $\dot{\bb{w}}$ given by \eqref{eq.w_adaptation_law} is affine in $\bb{u}$ as required in the statement of Problem \ref{prob.adaptation}, but, however, that because $\dot{\bb{w}}$ appears explicitly in \eqref{eq.b_ccbf_condition_unfiltered} (and therefore in \eqref{eq.Phi_func}) the adaptation law given by \eqref{eq.w_adaptation_law} with $b_{2c+1} = \beta$ requires solving a partial differential equation (PDE). Rather than attempt to solve this (possibly intractable) PDE, we break the algebraic loop by introducing $\bb{\mu}^f(t)$ and $\bb{\nu}^f(t)$ as filtered versions of $\bb{\mu}(t) \coloneqq \mu(t, \bb{w}(t), \bb{x}(t))$ and $\bb{\nu}(t) \coloneqq \nu(t, \bb{w}(t), \bb{x}(t))$. We assume that these filtered variables track their unfiltered counterparts sufficiently closely, as stated in the following.
\begin{Assumption}\label{ass.filter_errors}
    There exist bounded, positive constants $\eta_{\mu}$ and $\eta_{\nu}$, i.e., $0 < \eta_\mu,\eta_\nu < \infty$, such that the following relations hold:
    \begin{equation}
        \min_{t \in \mathcal{T}}\left[\frac{\partial H}{\partial \bb{w}}\left(\bb{\mu}^f(t) - \bb{\mu}(t)\right)\right] \geq -\eta_\mu, \label{eq.eta_mu}
    \end{equation}
    which implies that the deviation between the effect of the filtered signal $\bb{\mu}^f$ and true signal $\bb{\mu}$ on the C-CBF dynamics is bounded by $\eta_\mu$, and
    \begin{align}
        &\min_{t \in \mathcal{T}}\bigg[\sup_{\bb{u} \in \mathcal{U}}\left(\frac{\partial H}{\partial \bb{x}}g(\bb{x}(t)) + \frac{\partial H}{\partial \bb{w}}\bb{\nu}^f(t)\right)\bb{u} \label{eq.eta_nu} \\
        &\quad\quad\quad\quad\quad\; - \sup_{\bb{u} \in \mathcal{U}}\left(\frac{\partial H}{\partial \bb{x}}g(\bb{x}(t)) + \frac{\partial H}{\partial \bb{w}}\bb{\nu}(t))\right)\bb{u}\bigg] \geq -\eta_\nu, \nonumber
    \end{align}
    which similarly implies that the deviation (considering the worst input) between the effect of the filtered signal $\bb{\nu}^f$ and true signal $\bb{\nu}$ on the C-CBF dynamics is bounded by $\eta_\nu$.
\end{Assumption}

The above implicitly define worst-case error bounds $\eta_\mu$ and $\eta_\nu$ on how the filtered signals $\bb{\mu}^f$ and $\bb{\nu}^f$ affect $\dot H$ versus the true signals $\bb{\mu}$ and $\bb{\nu}$, and thus the inequality in \eqref{eq.b_ccbf_condition_unfiltered} may be satisfied by using $\bb{\mu}^f$, $\bb{\nu}^f$, $\eta_\mu$, and $\eta_\nu$ instead of $\bb{\mu}$ and $\bb{\nu}$. In practice, \eqref{eq.eta_mu} and \eqref{eq.eta_nu} may be satisfied via appropriate filter design.
We are now ready to define the final constraint function as
\begin{equation}\label{eq.b_ccbf_condition}
    b_{2c+1}(t, \bb{w}, \bb{x}) =\delta(t, \bb{w}, \bb{x}) - \left|q(t, \bb{w}, \bb{x})\right|^T\bb{u}_{max},
\end{equation}
which is directly \eqref{eq.b_ccbf_condition_unfiltered} with the filter variables (and error bounds), specifically
\begin{equation}\label{eq.delta}
\begin{aligned}
    \delta(t, \bb{w}, \bb{x}) = \eta_\mu + \eta_\nu -&\frac{\partial H}{\partial t} - \frac{\partial H}{\partial \bb{x}}f(\bb{x}) \\ 
    &- \frac{\partial H}{\partial \bb{w}}\bb{\mu}^f(t) - \alpha(H),
\end{aligned}
\end{equation}
and
\begin{equation}\label{eq.v_vector}
    q(t, \bb{w}, \bb{x}) \coloneqq \bb{L}_g^Tp_h(\bb{w}, \bb{x}) + \bb{\nu}^f(t)^Tp_w(\bb{w}, \bb{x}),
\end{equation}
with the C-CBF control matrix $\bb{L}_g$ (omitting the dependence on $\bb{x}$) defined by
\begin{align}
    \bb{L}_g &= \left[\frac{\partial h_1}{\partial \bb{x}} \cdots \frac{\partial h_c}{\partial \bb{x}}\right]^Tg(\bb{x}) \in \R^{c \times m}, \label{eq.Lg_matrix}
\end{align}
and the functions $p_h: \R^c \times \R^n \rightarrow \R^c$ and $p_w: \R^c \times \R^n \rightarrow \R^c$ given by
\begin{align}
    p_h(\bb{w}, \bb{x}) &= \left[\frac{\partial \phi(w_1,\bb{x})}{\partial h_1} \; \cdots \; \frac{\partial \phi(w_c,\bb{x})}{\partial h_c}\right]^T \label{eq.ph_vector} \\
    p_w(\bb{w}, \bb{x}) &= \left[\frac{\partial \phi(w_1,\bb{x})}{\partial w_1} \; \cdots \; \frac{\partial \phi(w_c,\bb{x})}{\partial w_c}\right]^T. \label{eq.pw_vector}
\end{align}
These modifications allow us to include $b_{2c+1}$ directly as a function constraining the optimal solution trajectory given by \eqref{eq.w_solution_trajectory} without needing to solve a complicated PDE.

% Before introducing our main result, we justify the following assumption on the rank of the C-CBF control matrix.
% \begin{Assumption}\label{ass.LGH_rank}
%     The matrix $\bb{L}_g$ given by \eqref{eq.Lg_matrix} has row-rank of at least one at all times, i.e $\rho(\bb{L}_g) \geq 1$ for all $t \geq 0$.
% \end{Assumption}

% In order for the row-rank $\rho(\bb{L}_g) \geq 1$, it must hold that $\frac{\partial h_i}{\partial \bb{x}}g(\bb{x}) \neq \mathbf{0}_{1 \times m}$ for some $i \in [c]$. For this to be satisfied, at least one constituent constraint function must have relative-degree one with respect to the system \eqref{eq.single_agent_system} (which is often reasonable, e.g., a speed constraint on a vehicle with acceleration control). Notably, any other constituent constraint function may have higher relative-degree. Without Assumption \ref{ass.LGH_rank}, i.e., if all constraint functions had relative-degree $r > 1$, control inputs would not affect the dynamics of any constraint function $h_i$.

We are now ready to state our main result.
\begin{Theorem}\label{thm.safe_adaptation}
    Consider $c\geq 1$ twice continuously differentiable constraint functions $h_i$ defining sets $\mathcal{S}_i(t)$ as in \eqref{eq.safe_set}, $\forall i \in [c]$, and the associated C-CBF candidate given by \eqref{eq.consolidated_cbf} with constituent weights $\bb{w}$. Suppose that Assumption \ref{ass.filter_errors} holds, and that $\mathrm{Int}(\mathcal{W}(t)) \neq \emptyset$ for all $t \geq 0$ with $\bb{w}(0) \in \mathrm{Int}(\mathcal{W}(0))$. Then, under the adaptation law \eqref{eq.w_adaptation_law}, the C-CBF candidate $H$ given by \eqref{eq.consolidated_cbf} is rendered valid, i.e., the condition \eqref{eq.c-cbf-condition} is satisfied for all $\bb{x}(t) \in \mathcal{D}(t)$ and $\delta \geq 0$, for all $t \geq 0$.
\end{Theorem}
\begin{proof}
    We will first derive an expression for $\dot{H}$, and then show that when $\dot{\bb{w}}$ is given by \eqref{eq.w_adaptation_law} there always exists a control input $\bb{u} \in \mathcal{U}$ such that $\sup_{\bb{u} \in \mathcal{U}}\dot{H}(t, \bb{w}, \dot{\bb{w}}, \bb{x}, \bb{u}) \geq -\alpha(H)$.
    
    First, observe that according to \eqref{eq.consolidated_cbf} the C-CBF candidate time-derivative $\dot{H}$ takes the following form:
    \begin{align}\nonumber
        \dot H &= -\sum_{i=1}^c\left(\frac{\partial \phi}{\partial h_i}\dot{h}_i + \frac{\partial \phi}{\partial w_i}\dot{w}_i\right) \nonumber \\
        &= -\left(p_h^T\dot{\bb{h}} + p_w^T\dot{\bb{w}}\right) \nonumber \\
        &= -\left(p_h^T\left(\bb{L}_t + \bb{L}_f + \bb{L}_g\bb{u}\right) + p_w^T(\bb{\mu} + \bb{\nu}\bb{u})\right), \nonumber
    \end{align}
    where $\bb{\mu} = \mu(t, \bb{w}, \bb{x})$ and $\bb{\nu} = \nu(t, \bb{w}, \bb{x})$ are given by \eqref{eq.ka_dot} and \eqref{eq.kb_dot} respectively, $\bb{L}_g$ by \eqref{eq.Lg_matrix}, $p_h$ and $p_w$ by \eqref{eq.ph_vector} and \eqref{eq.pw_vector}, and where (omitting the dependence on $\bb{x}$ and $\bb{w}$)
    \begin{align}
        \bb{L}_t &= \left[\frac{\partial h_1}{\partial t} \cdots \frac{\partial h_c}{\partial t}\right]^T \in \R^c, \nonumber \\
        \bb{L}_f &= \left[\frac{\partial h_1}{\partial \bb{x}} \cdots \frac{\partial h_1}{\partial \bb{x}}\right]^Tf(\bb{x}) \in \R^{c}. \nonumber 
    \end{align}
    Given $\alpha \in \mathcal{K}_{\infty}$, it follows that $\dot{H} + \alpha(H) = a + \bb{b}^T\bb{u}$, with $a = \alpha(H) - p_h^T(\bb{L}_t + \bb{L}_f) - p_w^T\bb{\mu}$ and $\bb{b} = -(\bb{L}_g^Tp_h + \bb{\nu}^Tp_w)$. Observe that given \eqref{eq.single_agent_system}  the adaptation law \eqref{eq.w_adaptation_law} may be expressed as
    \begin{equation*}\small
        \dot{\bb{w}} = -\nabla_{\bb{ww}}^{-1}\Phi\left(\bb{P}\nabla_{\bb{w}}\Phi + \nabla_{\bb{wx}}\Phi\dot{\bb{x}} + \nabla_{\bb{w}t}\Phi\right),
    \end{equation*}
    which, for $\Phi$ defined by \eqref{eq.Phi_func} defines a dynamical system of the form \eqref{eq.zdot_ode_optimization}. With $\Phi$ strongly, uniformly convex, $\mathrm{Int}(\mathcal{W}(t)) \neq \emptyset$ for all $t \geq 0$, and $\bb{w}(0) \in \mathrm{Int}(\mathcal{W}(0))$, it follows from Lemma \ref{lem.feasible_optimization_solution} that the solution remains within the feasible set for all time, i.e., $\bb{w}(t) \in \mathcal{W}(t)$, $\forall t \geq 0$, which implies that $b_j(t, \bb{w}(t), \bb{x}(t)) \leq 0$ for all $j \in [2c+1]$, $t \geq 0$. While constraints $j \in [c]$ encode that each $w_j > w_{min} > 0$ and constraints $j+c$ for $j \in [c]$ encode that each $w_j < w_{max} < \infty$, it may be seen by substituting \eqref{eq.delta} and \eqref{eq.v_vector} into \eqref{eq.b_ccbf_condition} that the final constraint function $b_{2c+1}$ encodes that 
    \begin{align*}
        \frac{\partial H}{\partial t} + \frac{\partial H}{\partial \bb{x}}f(\bb{x}) + \frac{\partial H}{\partial \bb{w}}\bb{\mu}^f &- \eta_\mu + \alpha(H) + \\ \bigg|\frac{\partial H}{\partial \bb{x}}g(\bb{x}) &+ \frac{\partial H}{\partial \bb{w}}\bb{\nu}^f\bigg|\bb{u}_{max} - \eta_\nu \geq 0.
    \end{align*}
    Then, by observing from the input constraint set \eqref{eq.input_constraint_set} that
    \begin{equation*}
        \bigg|\frac{\partial H}{\partial \bb{x}}g(\bb{x}) + \frac{\partial H}{\partial \bb{w}}\bb{\nu}^f\bigg|\bb{u}_{max} = \sup_{\bb{u} \in \mathcal{U}}\left(\frac{\partial H}{\partial \bb{x}}g(\bb{x}) + \frac{\partial H}{\partial \bb{w}}\bb{\nu}^f\right)\bb{u}
    \end{equation*}
    and taking by \eqref{eq.eta_mu} and \eqref{eq.eta_nu} from Assumption \ref{ass.filter_errors}, it follows that that $b_{2c+1} \leq 0$ implies that
    \begin{align*}
        \dot H = \frac{\partial H}{\partial t} + \frac{\partial H}{\partial \bb{x}}f(\bb{x}) &+ \frac{\partial H}{\partial \bb{w}}\mu(t, \bb{w}, \bb{x}) + \alpha(H) \\ +&\sup_{\bb{u} \in \mathcal{U}}\bigg[\bigg(\frac{\partial H}{\partial \bb{x}}g(\bb{x}) + \frac{\partial H}{\partial \bb{w}}\nu(t, \bb{w}, \bb{x})\bigg)\bb{u}\bigg], \\
        = \frac{\partial H}{\partial t} + \frac{\partial H}{\partial \bb{x}}f(\bb{x}) &+ \frac{\partial H}{\partial \bb{w}}\dot{\bb{w}} +  \sup_{\bb{u} \in \mathcal{U}}\frac{\partial H}{\partial \bb{x}}g(\bb{x})\bb{u} \geq -\alpha(H),
    \end{align*}
    and therefore there always exists a control input $\bb{u} \in \mathcal{U}$ such that the C-CBF condition \eqref{eq.c-cbf-condition} is satisfied. Thus, the adaptation law \eqref{eq.w_adaptation_law} renders the C-CBF candidate valid. This completes the proof.
\end{proof}

The above result is predicated on the following assumptions holding. First, the convexity condition imposed on $\Phi$ as outlined in Assumption \ref{ass.psi_convex}: this may require monitoring the eigenvalues of $\nabla_{\bb{ww}}\Phi$ online and increasing $s$ if necessary, but in our experiments we have found this to be unproblematic. Second, and most challenging, the filter design for $\bb{\mu}(t)$ and $\bb{\nu}(t)$ must be sufficiently accurate as detailed in Assumption \ref{ass.filter_errors}. The robustness margins $\eta_\mu$ and $\eta_\nu$ may be difficult to determine a priori without knowledge of how the system will evolve over time. In our numerical experiments, we used a trial and error process with first-order low-pass filters and determined that $\eta_\mu,\eta_\nu \leq 0.05$ typically worked well for the case studies undertaken. The third and final set of assumptions involves the feasible region $\mathcal{W}$, namely that $\mathrm{Int}(\mathcal{W}(t)) \neq \emptyset$ and that $\bb{w}(0) \in \mathcal{W}(0)$. Ensuring that the initial weights belong to the feasible set is fairly straightforward: we choose a set of weights and if they do not belong to the feasible set we adapt them using only the correction term of the adaptation law $\dot{\bb{w}}_0 = -\nabla_{\bb{ww}}^{-1}\Phi\bb{P}\nabla_{\bb{w}}\Phi$ until $\bb{w}_0 \in \mathrm{Int}(\mathcal{W}(0))$. Ensuring that $\mathrm{Int}(\mathcal{W}(t)) \neq \emptyset$ for all $t \in \mathcal{T}$ is currently an open problem, and thus we must assume that it is true. What our result in Theorem \ref{thm.safe_adaptation} does imply, however, is that if the feasible region has an interior point then the adaptation law \eqref{eq.w_adaptation_law} will find it. If no interior feasible point exists, then we may use various relaxation strategies on the constraint functions $b_j$, $j \in [2c]$ (which excludes the constraint function encoding the C-CBF condition) in order to consider an enlarged feasible region and ensure that our adaptation law is well-defined. We refer the interested reader to e.g., the discussion in \cite[Section III.C]{Fazlyab2018Prediction}.

\begin{Remark}
    Notably absent from the set of requisite assumptions is any condition on the relative-degree of the state constraint functions $h_i$. For our method it is actually acceptable for the C-CBF control matrix $\bb{L}_g$ to consist only of zeros because the parameter adaptation law $\dot{\bb{w}}$ contributes to the satisfiability of the C-CBF condition \eqref{eq.c-cbf-condition} via the term $\frac{\partial H}{\partial \bb{w}}\nu(t, \bb{w}, \bb{x})\bb{u}$. In this case, the adaptation of the weights would be fully responsible for ensuring that \eqref{eq.c-cbf-condition} holds. This does not come without tradeoffs, however, since it may be more difficult to satisfy the assumption on the initial condition that $\bb{w}(0) \in \mathcal{W}(0)$ (which may restrict the set of allowable initial states) or to ensure that $\mathrm{Int}(\mathcal{W}(t)) \neq \emptyset$.
\end{Remark}

\section{Consolidated CBF based Control Synthesis}\label{sec.ccbf_control}

In this section, we formally introduce a controller based on the proposed consolidated control barrier function and discuss its present limitations. 

\subsection{C-CBF-QP Control Design}
We previously reviewed the CBF-QP control law given by \eqref{eq.cbf_qp_controller} and suggested that our C-CBF could be directly inserted into such a framework for certified constrained control design. The C-CBF-QP controller can be described by
\begin{equation}\label{eq.u_solution_trajectory}
\begin{aligned}
    \bb{u}^*(t) &= \argmin_{\bb{u} \in \R^m} \; J_u(t, \bb{x}, \bb{u}), \;\; \textrm{s.t.} \;\; \bb{u} \in \mathcal{U}_H(t, \alpha),
\end{aligned}
\end{equation}
where the objective function $J_u: \mathcal{T} \times \R^n \times \R^m \rightarrow \R$ is given by
\begin{equation*}
    J_u(t, \bb{x}, \bb{u}) = \frac{1}{2}\|\bb{u} - \bb{u}_0(t, \bb{x})\|^2,
\end{equation*}
and the feasible region $\mathcal{U}_H$ is the C-CBF control set given by \eqref{eq.ccbf_control_set}, defined more concisely here as
\begin{equation*}
    \mathcal{U}_H(t,\alpha) = \{\bb{u} \in \R^m \mid d_k(t, \bb{w}, \bb{x}, \bb{u}) \leq 0, \forall k \in [2m+1]\}
\end{equation*}
for input constraint functions
\begin{align}
    d_k(t, \bb{w}, \bb{x}, \bb{u}) &= \bar{u}_k - u_k, \; \forall k \in \{1,\hdots,m\}, \label{eq.d_umin} \\
    d_k(t, \bb{w}, \bb{x}, \bb{u}) &= u_k - \bar{u}_k, \; \forall k \in \{m+1,\hdots,2m\},\label{eq.d_umax}
\end{align}
and the C-CBF constraint function
\begin{equation*}
\begin{aligned}
        d_{2m+1}(t, \bb{w}, \bb{x}, \bb{u}) = -\alpha(H) -& \frac{\partial H}{\partial t} - \frac{\partial H}{\partial \bb{x}}f(\bb{x}) -\frac{\partial H}{\partial \bb{w}}\dot{\bb{w}} \\
        & -\frac{\partial H}{\partial \bb{x}}g(\bb{x})\bb{u}.
\end{aligned}
\end{equation*}
In the following result, we highlight how this proposed controller results in complete constraint satisfaction when $\dot{\bb{w}}$ is given by \eqref{eq.w_adaptation_law}.
\begin{Theorem}\label{thm.safe_control}
    Consider the set of $c \geq 1$ constraint functions $\mathcal{H}$ given by \eqref{eq.constraint_function_set} and the complete constraint set $\mathcal{S}(t)$ given by \eqref{eq.complete_constraint_set}. Assume that the premises of Theorem \ref{thm.safe_adaptation} hold, and that $\dot{\bb{w}}$ is given by \eqref{eq.w_adaptation_law}. Then, under the control law \eqref{eq.u_solution_trajectory} the state remains inside the complete constraint set and thus all constraints encoded via functions $h \in \mathcal{H}$ are satisfied for all time, i.e., $\bb{x}(t) \in \mathcal{S}(t)$, $\forall t \in \mathcal{T}$ under control policy $\bb{u}^*(t)$ given by \eqref{eq.u_solution_trajectory}.
\end{Theorem}
\begin{proof}
    First, we have from Theorem \ref{thm.safe_adaptation} and its requisite assumptions that when $\dot{\bb{w}}$ is given by \eqref{eq.w_adaptation_law} then the C-CBF control set is non-empty at all times, i.e., $\mathcal{U}_H(t, \alpha) \neq \emptyset$, $\forall t \in \mathcal{T}$. This implies that the optimization problem \eqref{eq.u_solution_trajectory} is feasible at all times, from which it follows that $\bb{u}^*(t) \in \mathcal{U}_H(t,\alpha)$, $\forall t \in \mathcal{T}$. Thus, we have that the function $H$ is a valid CBF with respect to the set $\mathcal{D}(t)$ for the system \eqref{eq.single_agent_system}. It then follows from Theorem \ref{thm.cbf_invariance} that $\bb{x}(t) \in \mathcal{D}(t)$, $\forall t \in \mathcal{T}$. Since the set $\mathcal{D}(t)$ is a subset of the complete constraint set $\mathcal{S}(t)$ at all times, i.e., $\mathcal{D}(t) \subset \mathcal{S}(t)$, $\forall t \in \mathcal{T}$, we have that $\bb{x}(t) \in \mathcal{S}(t)$ and therefore $h(t, \bb{x}(t)) \geq 0$, $\forall h \in \mathcal{H}$, for all $t \in \mathcal{T}$.
\end{proof}

\subsection{Discussion}
From the above result, it is evident that the purpose behind the parameter adaptation law $\dot{\bb{w}}$ given by \eqref{eq.w_adaptation_law} is to enable the feasibility of the QP controller described by \eqref{eq.u_solution_trajectory}. In this sense, we have transferred the feasibility problem from one optimization problem (the control law \eqref{eq.u_solution_trajectory}) to another (\eqref{eq.w_solution_trajectory}, the basis for the adaptation law \eqref{eq.w_adaptation_law}) without providing any formal feasibility guarantees. There is some precedent for this type of approach, e.g., in \cite{xiao2022sufficient} the authors augment the QP-control law with an additional constraint encoding feasibility under input constraints and (at some point) assume feasibility of this new problem. In our case, we assume that the foundational optimization problem for our proposed adaptation law satisfies Slater's condition, i.e., that the feasible region contains an interior point. Without this assumption, analyzing the stability properties of the interconnected dynamical system
\begin{equation*}
    \dot{\bb{\xi}} = \begin{bmatrix}
        \dot{\bb{x}} \\ \dot{\bb{w}}
    \end{bmatrix} = \begin{bmatrix}
        f(\bb{x}) \\ \mu(t, \bb{w}, \bb{x})
    \end{bmatrix} + \begin{bmatrix}
        g(\bb{x}) \\ \nu(t, \bb{w}, \bb{x})
    \end{bmatrix}\bb{u}^*,
\end{equation*}
is a difficult task. We recognize that this is a limitation to the proposed method, although we see potential in this line of research in that interior point methods similar to that used to derive the adaptation law \eqref{eq.w_adaptation_law} from optimization problem \eqref{eq.w_solution_trajectory} may be used to introduce a new dynamical system $\dot{\bb{u}} = \upsilon(t,\bb{w},\bb{x},\bb{u})$ for the input $\bb{u}$ to systems \eqref{eq.single_agent_system} and \eqref{eq.w_adaptation_law}, as outlined in what follows.

\begin{figure}[!ht]
    \centering
        \includegraphics[clip,width=1\linewidth]{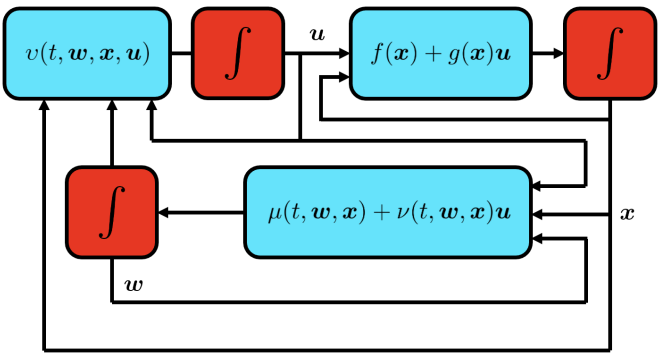}
    \caption{Block diagram for the $\dot{\bb{\chi}}$ dynamics described by \eqref{eq.chi_system_dynamics}.}\label{fig.block_diagram}
\end{figure}

Given that the objective function $J_u$ for \eqref{eq.u_solution_trajectory} is twice continuously differentiable and strongly, uniformly convex in $\bb{u}$ and that the constraint functions $d_k$, $k \in [2m+1]$ are convex (because they are affine) in $\bb{u}$, we may again use log-barriers to introduce the approximate optimal control trajectory as
\begin{equation}\label{eq.approx_u_solution_trajectory}
\begin{aligned}
    \hat{\bb{u}}^*(t) &= \argmin_{\bb{u} \in \R^m} \; \Omega(t, \bb{w}, \bb{x}, \bb{u}),
\end{aligned}
\end{equation}
where
\begin{equation}
    \Omega(t, \bb{w}, \bb{x}, \bb{u}) \coloneqq J_u(t, \bb{x}, \bb{u}) - \frac{1}{s}\sum_{k=1}^{2m+1}\log\big(-d_k(t, \bb{w}, \bb{x}, \bb{u})\big),
\end{equation}
for barrier parameter $s>0$. Assume (as for $\Psi$ in \eqref{eq.approximate_solution}, $\Phi$ in \eqref{eq.what_solution_trajectory}) that Assumption \ref{ass.psi_convex} holds for $\Omega$, i.e., that there exists $s>0$ that renders $\Omega$ strongly, uniformly convex in $\bb{u}$, and consider the following dynamical system
\begin{equation}\label{eq.u_dot_system}
\begin{aligned}
    \dot{\bb{u}} &= \upsilon(t, \bb{w}, \bb{x}, \bb{u}), \\
    &\coloneqq -\nabla_{\bb{uu}}^{-1}\Omega\bigg(\bb{B}\nabla_{\bb{u}}\Omega + \nabla_{\bb{uw}}\Omega\dot{\bb{w}}
        + \nabla_{\bb{ux}}\Omega\dot{\bb{x}} + \nabla_{\bb{u}t}\Omega\bigg), \\
    \bb{u}&(0) \in \mathcal{U}_H(0, \alpha),
\end{aligned}
\end{equation}
with $\bb{B} \in \R^{m\times m}$ a positive-definite gain matrix affecting the strength of the correction term, $\dot{\bb{w}}$ the adaptation law given by \eqref{eq.w_adaptation_law}, and $\dot{\bb{x}}$ the state dynamics given by \eqref{eq.single_agent_system}. In this case, the control law could then be given as the solution $\bb{u}(t)$ to the system \eqref{eq.u_dot_system}, i.e.,
\begin{equation}\label{eq.proposed_controller}
    \bb{u}(t) = \varphi_t(\bb{u}(0)),
\end{equation}
where $\varphi^t$ is the flow map of \eqref{eq.u_dot_system} such that $t \mapsto \varphi^t(\bb{u}(0))$ solves \eqref{eq.u_dot_system} for initial condition $\bb{u}(0)$. 
% In practice, we determine $\dot{\bb{u}}$ pointwise in time and solve for $\bb{u}(t)$ using finite-difference methods.
% \begin{Remark}
%     It is shown that the solution to \eqref{eq.u_dot_system} converges asymptotically only to the \texit{approximate} optimal control trajectory $\hat{\bb{u}}^*(t)$ rather than the true optimal solution trajectory $\bb{u}^*(t)$. This does not necessarily result in decreased system performance, however, as generally speaking the class of CBF-QP control laws is only optimal pointwise-in-time rather than over a time interval. There is no guarantee of optimality over an operation period, and therefore our proposed controller does not sacrifice optimality over any specific time interval with respect to the standard class of CBF-QP control laws.
% \end{Remark}
% The only additional assumption required for Theorem \ref{thm.safe_control} is the convexity of the augmented cost function $\Omega$. The viability of a control input satisfying the C-CBF condition under input constraints is a direct result of the adaptation law \eqref{eq.w_adaptation_law}, which means that a solution trajectory \eqref{eq.u_solution_trajectory} is guaranteed to exist. 
The motivation for using interior point methods to introduce the dynamical system \eqref{eq.u_dot_system} is that it enables the expression of the state $\bb{x}$, weight $\bb{w}$, and control $\bb{u}$ dynamics as one interconnected dynamical system:
\begin{equation}\label{eq.chi_system_dynamics}
    \dot{\bb{\chi}} = F(t, \bb{\chi}) = \begin{bmatrix}
        \dot{\bb{x}} \\ \dot{\bb{w}} \\ \dot{\bb{u}}
    \end{bmatrix} = \begin{bmatrix}
        f(\bb{x}) + g(\bb{x})\bb{u} \\ \mu(t, \bb{w}, \bb{x}) + \nu(t, \bb{w}, \bb{x})\bb{u} \\ \upsilon(t, \bb{w}, \bb{x}, \bb{u})
    \end{bmatrix}
\end{equation}
where $\bb{\chi} = (\bb{x}, \bb{w}, \bb{u}) \in \R^{n+c+m}$ is the state, and for which a block diagram is provided in Figure \ref{fig.block_diagram}. Though this is not the first time that a continuous-time interior point method has been used in place of a CBF-QP control law (see e.g., \cite{wang2022suboptimal}), we are unaware of any work which uses such an approach for joint adaptation and control for constrained control design. This permits the use of available tools like Barbalat's lemma \cite[Lem. 8.2]{khalil2002nonlinear} for stability analysis of time-varying systems, and we hope to investigate these properties in the future.

%%%%%%%%%%%%%%%%%%%%%%%%%%%%%%%%%%%%%%%%%%%%%%%%%%%%%%%%%%
%************ Numerical Studies **************%
%%%%%%%%%%%%%%%%%%%%%%%%%%%%%%%%%%%%%%%%%%%%%%%%%%%%%%%%%%
\section{Numerical Case Studies}\label{sec.case_studies}
In this section, we highlight the performance of our proposed adaptive C-CBF controller on two examples: an illustrative one-dimensional problem, and a reach-avoid scenario inspired by mobile robots with which we compare our method to several baseline controllers.

\subsection{One-Dimensional Nonlinear System}
For demonstration purposes we consider the following nonlinear, control-affine dynamical system:
\begin{equation}
    \dot x = x(e^{{kx}^2}-1) + (4 - x^2)u, \quad x(0) = 0,
\end{equation}
where $x \in \R$ denotes the state, $u \in \mathcal{U} \subset \R$ the control input, and where $k = \ln(2) / p^2$ with $p=1.5616$. It may be seen that the uncontrolled system has a unique (unstable) equilibrium at $x=0$, and that there are two distinct control singularities at $x=2$ and $x=-2$ respectively in the sense that the control $u$ has no effect on the dynamics at these points. The control objective is to track a known trajectory 
\begin{equation*}
    x^*(t) = 4\sin(2\pi \frac{t}{5} + \theta),
\end{equation*}
where we consider separately $\theta=0$ and $\theta=\pi$ (i.e., the sinusoid evolving in both positive and negative directions) subject to the state constraints $x \leq 2$ and $x \geq -2$, which may be encoded via constraint functions
\begin{align*}
    h_1(x) &= 2 - x, \\
    h_2(x) &= x + 2.
\end{align*}
It is evident that perfect tracking of the desired state would result in violations of both of the above constraint functions as the signal $x^*(t)$ oscillates. We take the input constraint set to be $\mathcal{U} = [-1,1]$, and note that this creates two infinite potential wells in the following sense: if it occurs that $2 > x > p_1 = p$ (resp. $-2 < x < p_2 = -p$) then $x$ will escape to $\infty$ (resp. $-\infty$) and violate $h_1$ (resp. $h_2$) in doing so. In short, though the state constraints are $x \leq 2$ and $x \geq -2$, due to the input constraint set $\mathcal{U}$ these constraints are guaranteed to be violated if either $x > p$ or $x < -p$ regardless of the control applied thereafter. Our adaptive C-CBF controller has the ability to protect against these types of potential wells via online parameter adaptation, and we now show how it manages to continuously satisfy both constraints $h_1$ and $h_2$ simultaneously despite seeking to track an unsafe nominal trajectory.

We simulated the tracking problem using three different iterations of the C-CBF-QP controller \eqref{eq.u_solution_trajectory} with a proportional control law tracking the desired trajectory $x^*(t)$ as a nominal input. The controllers differed in their aggressiveness as encoded via the following class-$\mathcal{K}$ functions appearing in the C-CBF condition: $\alpha(H) = \gamma \cdot H^3$ with $\gamma = 0.01, 0.1, 1.0$. This made for a total of six trials: three controllers tested each on two separate desired trajectories. In all six scenarios the initial weights were $\bb{w}(0) = (2.38, 2.38)$, determined by adapting according to the correction term from $(1,1)$ as outlined in Section \ref{sec.consolidated_cbf}. For our C-CBF $H$ we used \eqref{eq.consolidated_cbf} with decaying exponentials, i.e., $\phi(r,s) = e^{-rs}$. For our adaptation law \eqref{eq.w_adaptation_law}, we used $w_{min} = 0.01$, $w_{max} = 50.0$, and chose $\bb{P} = 100\bb{I}$.

\begin{figure}
    \centering
        \includegraphics[clip,width=1\linewidth]{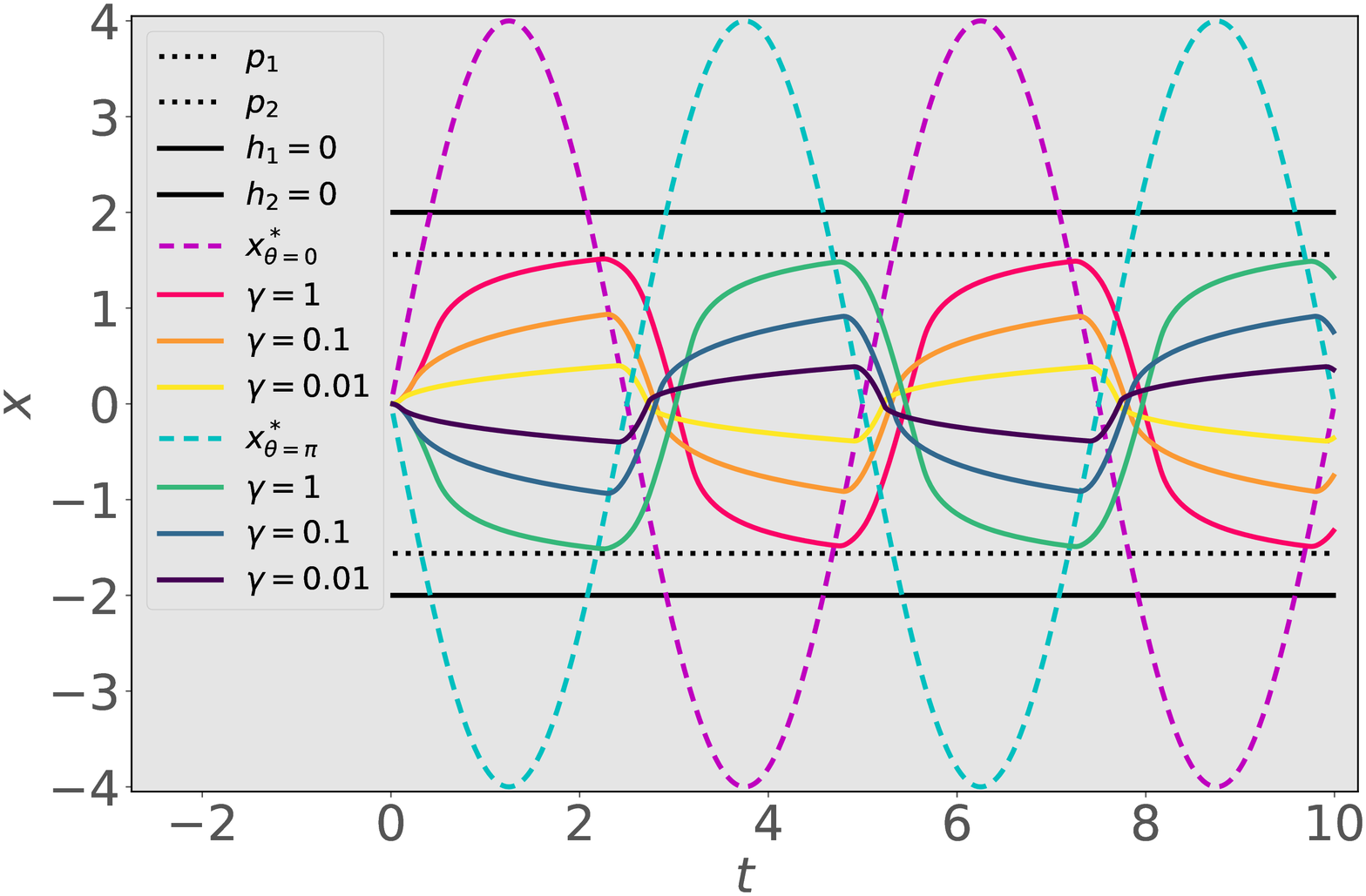}
    \caption{Evolution of the state $x$ for the 1-D trajectory tracking problem using the C-CBF controller with 3 different class $\mathcal{K}_\infty$ functions of the form $\alpha(H) = \gamma \cdot H^3$ for the given $\gamma$ values.}\label{fig.nonlinear1d_state}
\end{figure}

\begin{figure}
    \centering
        \includegraphics[clip,width=1\linewidth]{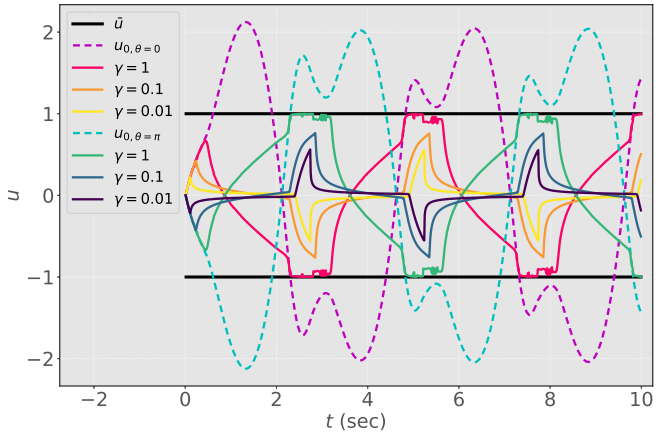}
    \caption{Control inputs using the C-CBF controller with 3 different class $\mathcal{K}_\infty$ functions for the 1-D trajectory tracking problem.}\label{fig.nonlinear1d_controls}
\end{figure}

\begin{figure}
    \centering
        \includegraphics[clip,width=1\linewidth]{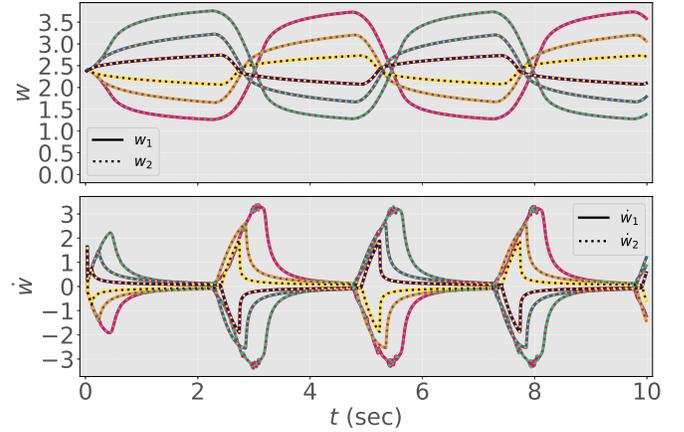}
    \caption{Constituent constraint function weight $\bb{w}$ and weight derivative $\dot{\bb{w}}$ trajectories for the 1-D trajectory tracking problem. Colors identical to those found in legend of Figure \ref{fig.nonlinear1d_ccbf}.}\label{fig.nonlinear1d_parameters}
\end{figure}

\begin{figure}
    \centering
        \includegraphics[clip,width=1\linewidth]{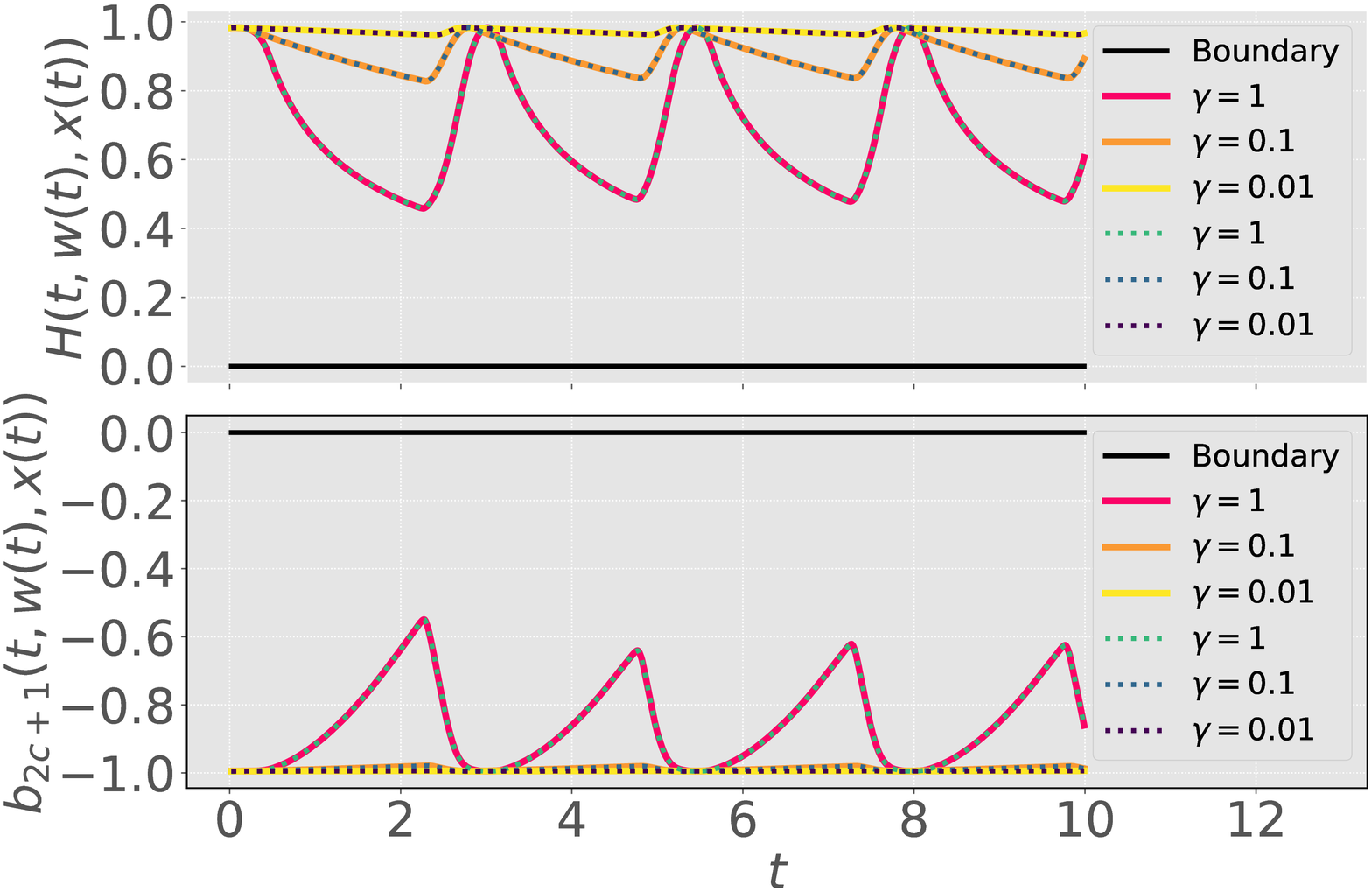}
    \caption{Values of the C-CBF $H$ and the constraint function $b_{2c+1}$ for C-CBF condition satisfiability for the 6 different C-CBF controlled trajectory tracking scenarios.}\label{fig.nonlinear1d_ccbf}
\end{figure}

The evolution of the state in each trial is shown in Figure \ref{fig.nonlinear1d_state}, from which it may be seen that not only are the position constraints satisfied at all times despite the unsafe reference path but also every iteration of the C-CBF controller succeeds in preventing the system from entering an infinite potential wells. As can be expected, the controllers equipped with more aggressive parameters $\gamma$ approached the boundary of the potential wells more closely. It is evident from Figure \ref{fig.nonlinear1d_controls} that the proposed controller modifies the nominally unsafe inputs in advance of any danger. Note that while we do not provide any guarantees of Lipschitz continuous control inputs or weights, the approach produced both inputs and weights that varied smoothly, as seen also in Figure \ref{fig.nonlinear1d_parameters}. It is worth highlighting the symmetry of the adaptation for identical controllers tracking an inverted reference trajectory, which makes sense given the symmetry of the dynamics about the equilibrium point and the identical initial conditions. Figure \ref{fig.nonlinear1d_ccbf} verifies that the C-CBF $H$ remained non-negative at all times and that the C-CBF condition remained satisfiable at all times under the input constraints thanks to the proposed adaptation law \eqref{eq.w_adaptation_law}.


\subsection{Mobile Robot Reach-Avoid Scenario}

In our second numerical case study, we consider a (car-like) mobile robot with bicycle dynamics seeking to reach a target location within a prescribed time in the presence of static obstacles, a speed limit, and steering constraints.
Let $s_0$ be a local origin in an inertial frame $\mathcal{I}$, and suppose that the robot dynamics may be described by the following dynamic extension of the non-holonomic bicycle model from \cite[Ch. 2]{Rajamani2012VDC}:
\begin{subequations}\label{eq.dynamic_bicycle_model}
\begin{align}
    \dot{x} &= v\left(\cos{\psi} - \sin{\psi}\tan{\beta}\right), \label{eq: dyn x} \\
    \dot{y} &= v\left(\sin{\psi} + \cos{\psi}\tan{\beta}\right), \label{eq: dyn y} \\
    \dot{\psi} &= \frac{v}{l_r}\tan{\beta}, \label{eq: dyn psi} \\
    \dot{\beta} &= \omega, \\
    \dot{v} &= a,
\end{align}
\end{subequations}
where $x$ and $y$ denote the position (in m) of the center of gravity (c.g.) of the robot with respect to $s_0$, $\psi$ is the orientation (in rad) of its body-fixed frame, $\mathcal{B}$, with respect to $\mathcal{I}$, $\beta$ is the slip angle\footnote{$\beta$ is related to the steering angle $\delta$ via $\tan{\beta} = \frac{l_r}{l_r+l_f}\tan{\delta}$, where $l_f+l_r$ is the wheelbase with $l_f$ (resp. $l_r$) the distance from the c.g. to the center of the front (resp. rear) wheel.} (in rad) of the c.g. of the vehicle relative to $\mathcal{B}$ (with $|\beta|<\frac{\pi}{2}$), and $v$ is the velocity of the rear wheel with respect to $\mathcal{I}$. The state is denoted $\bb{z} = [x\; y \; \psi \; \beta \; v]^T$, and the control input is $\bb{u}=[\omega \; a]^T$, where $\omega$ is the angular velocity (in rad/s) of $\beta$ and $a$ is the acceleration of the rear wheel (in m/s$^2$). The objective of the robot is to reach a neighborhood of the goal location $(x_g, y_g) = (2, 2)$ from initial condition $\bb{z}_0 = [0, \; 0, \; \mathrm{arctan}\big((y_g - y) / (x_g - x)\big) \; 0.5, \; 0.0]^T$ while obeying the following 8 constraints: avoid five circular obstacles, obey the speed limit of $S=2$m/s, obey the slip angle limit $B = \pi/3$, and reach the goal set within $T=5$sec. The constraint functions encoding avoidance of the i$^{th}$ circular obstacle are given by
\begin{equation}
    h_i(t, \bb{z}) = (x - c_{x,i})^2 + (y - c_{y,i})^2 - R_i^2, \; \forall i \in [5], \nonumber
\end{equation}
for radius $R_i > 0$ and center point $(c_{x,i},c_{y,i})$. The speed and slip angle constraints are 
\begin{align*}
    h_6(t, \bb{z}) &= S^2 - v^2, \\
    h_7(t,\bb{z}) &=  B^2 - \beta^2,
\end{align*}
and the time-specification is defined by
\begin{equation*}
    h_8(t, \bb{z}) = R_g^2 + R_i^2\left(1 - \frac{t}{T}\right)^2 - (x - x_g)^2 - (y - y_g)^2,
\end{equation*}
for the goal set centered at $(x_g,y_g)$ with radius $R_g=0.1$m and shrinking radius $R_i=4$m. Note that with the exception of constraint functions $h_6$ and $h_7$, all other constraint functions have relative-degree two with respect to the dynamics \eqref{eq.dynamic_bicycle_model}, which means that the C-CBF control matrix $\bb{L}_g$ has six rows consisting only of zeros. If these functions were used as CBF candidates directly, their control terms would be zero at all times, i.e., $L_gh_i = \mathbf{0}$ for $i \in \{1, 2, 3, 4, 5, 8\}$, $\forall t \geq 0$. The C-CBF is then of the form \eqref{eq.consolidated_cbf}, where again use exponentials $\phi(r,s) = e^{-rs}$ with initial weights $\bb{w}(0) = 1.0 \cdot \mathbf{1}_{8 \times 1}$. We use the function $\alpha(H) = H$. 

\begin{figure}[!ht]
    \centering
        \includegraphics[clip,width=1\linewidth]{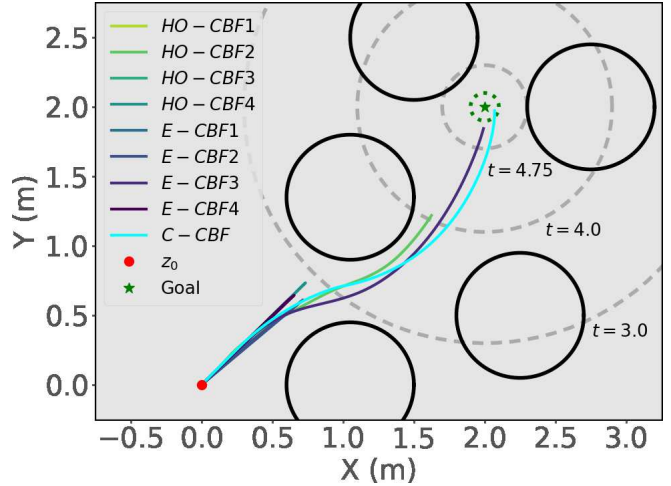}
    \caption{XY paths for bicycle robots seeking to reach the goal set (dotted green circle) within a prescribed time in the presence of static obstacles (black circles), a speed limit, and a slip angle constraint. Gray dashed lines indicate the reach constraint set at various times. All 6 HO- and E-CBF QP controllers become infeasible (though one E-CBF comes very close to the goal), while our adaptive C-CBF controller guarantees sufficient control authority for the feasibility of \eqref{eq.u_solution_trajectory} under input constraints and thus reaches the goal.}\label{fig.toy_example_paths}
\end{figure}
\begin{figure}[!ht]
    \centering
        \includegraphics[clip,width=1\linewidth]{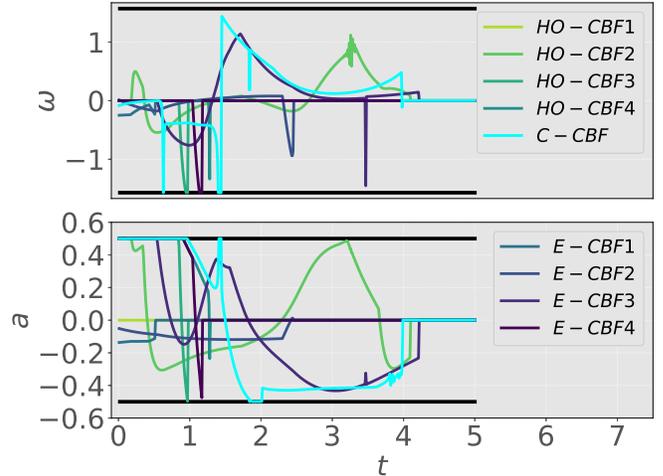}
    \caption{Control trajectories for the bicycle robot reach-avoid problem.}\label{fig.toy_example_controls}
\end{figure}
\begin{figure}[!ht]
    \centering
        \includegraphics[clip,width=1\linewidth]{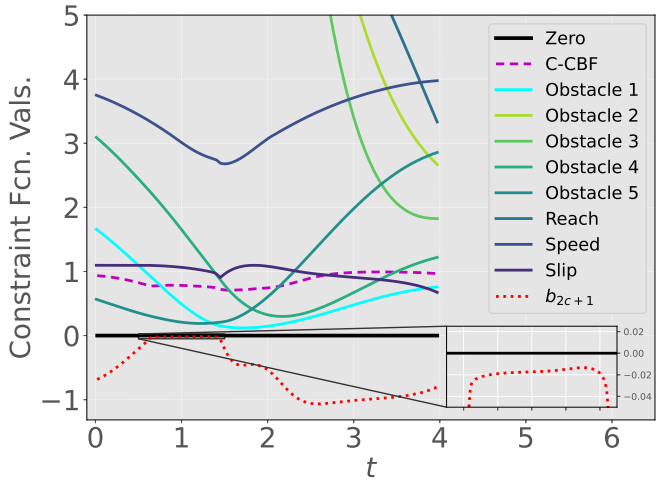}
    \caption{Values of the constraint functions $h_i$, $\forall i \in [8]$ (solid), the C-CBF $H$ (dashed), and constraint function $b_{2c+1}$ given by \eqref{eq.b_ccbf_condition} (dotted) with respect to time for the C-CBF controlled bicycle robot in the reach-avoid problem.}\label{fig.toy_example_cbfs}
\end{figure}
\begin{figure}[!ht]
    \centering
        \includegraphics[clip,width=1\linewidth]{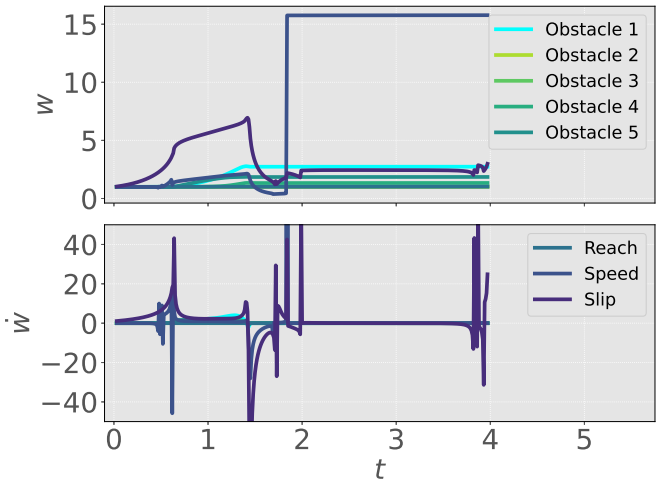}
    \caption{C-CBF weights $\bb{w}$ and their derivatives $\dot{\bb{w}}$ using adaptation law \eqref{eq.w_adaptation_law} for the bicycle robot reach-avoid problem.}\label{fig.toy_example_weights}
\end{figure}

For comparison against existing works, we simulate the following: a HO-CBF-QP controller (with guaranteed input constraint satisfaction for a single CBF) proposed by \cite{Breeden2021InputConstraints}, and an exponential (E-) CBF-QP controller (a subset of the class of HO-CBF controllers from \cite{Xiao2022HighOrder}) introduced in \cite{nguyen2016exponential}, both of which require reformulating constraint functions $h_i$ for $i \in \{1, 2, 3, 4, 5, 8\}$ as high-order CBFs. We simulate each class of controller over a time interval of $\mathcal{T}=[0,5]$ sec at a timestep of $\Delta t = 0.01$ sec and under four different class-$\mathcal{K}_\infty$ functions with varying levels of conservatism, from most (e.g., HO-CBF-1) to least conservative (e.g., HO-CBF-4). For the full list of parameters, as well as code and simulation videos, visit our Github repository\footnote{Github repo: \href{https://github.com/6lackmitchell/CCBF-Control}{https://github.com/6lackmitchell/CCBF-Control}}. The resulting paths and controls applied by the simulated bicycle robots are shown in Figures \ref{fig.toy_example_paths} and \ref{fig.toy_example_controls} respectively. Although the HO-CBF-QP controller proposed by \cite{Breeden2021InputConstraints} guarantees constraint adherence under input constraints for one constraint function, under multiple constraint functions the QP controller becomes infeasible in all four simulated trials (doing so, in fact, almost immediately for the most conservative HO-CBF-1 case) and therefore does not reach the goal. The E-CBF-QP controller has no guaranteed input constraint satisfaction and consequently takes more aggressive control actions. Notice, however, that the performance is highly sensitive to the choice of E-CBF gains: too conservative (E-CBF-1) and the QP quickly becomes infeasible; too aggressive (E-CBF-4) and the QP becomes infeasible when it no longer has sufficient control authority to avoid a collision with the first obstacle. It is clear that for this example both trial and error and expert knowledge are required to tune the E-CBF-QP controller to solve the problem. 

For our adaptive C-CBF controller, however, the assumption that the feasible region $\mathcal{W}(t)$ has an interior point at all times holds, and thus our proposed adaptation law finds a weighting of the constraint functions that allows the robot to safely reach the goal within the prescribed time. Figure \ref{fig.toy_example_cbfs} shows that all constraints were satisfied, and that feasibility of the C-CBF condition is preserved despite six of the eight constraints having relative-degree two with respect to the system dynamics (and thus contributing nothing to the $\frac{\partial H}{\partial \bb{x}}g(\bb{x})$ term in the C-CBF dynamics). Smoothly weighting the various constraint functions, as highlighted in Figure \ref{fig.toy_example_weights}, allows our controller to take smooth control actions without sudden switching or oscillatory behavior, and thus the adaptive C-CBF controlled vehicle reaches the goal while satisfying every constraint throughout the duration of the maneuver.

\section{Conclusion}\label{sec.conclusion}
In this paper, we addressed the problem of safe-control synthesis under multiple spatiotemporal and input constraints via a novel consolidated control barrier function. We first synthesized one candidate CBF from many weighted constituents, and then proposed a weight-adaptation law, inspired by interior-point methods, under which the C-CBF is rendered valid for constrained control design. We proposed the use of a C-CBF-QP control law, and proved its feasibility under the developed adaptation law. To validate the approach, we demonstrated its success on an unstable 1-D example and further highlighted its utility for generic constraints of arbitrary relative-degree on a robot-inspired reach-avoid problem subject to spatiotemporal constraints.

In the future, we plan to explore the stability properties of the complete, interconnected dynamical system given by \eqref{eq.chi_system_dynamics} in hopes of robustifying the constraint satisfaction guarantees. Additional future directions may also include the use of alternative interior point methods to improve the asymptotic convergence properties of the current approach, generalizing the class of input constraints considered, and exploring the optimality of the adaptation and control laws over time intervals.

\bibliographystyle{IEEEtran}
\bibliography{root}

\begin{thebibliography}{10}
\providecommand{\url}[1]{#1}
\csname url@samestyle\endcsname
\providecommand{\newblock}{\relax}
\providecommand{\bibinfo}[2]{#2}
\providecommand{\BIBentrySTDinterwordspacing}{\spaceskip=0pt\relax}
\providecommand{\BIBentryALTinterwordstretchfactor}{4}
\providecommand{\BIBentryALTinterwordspacing}{\spaceskip=\fontdimen2\font plus
\BIBentryALTinterwordstretchfactor\fontdimen3\font minus
  \fontdimen4\font\relax}
\providecommand{\BIBforeignlanguage}[2]{{%
\expandafter\ifx\csname l@#1\endcsname\relax
\typeout{** WARNING: IEEEtran.bst: No hyphenation pattern has been}%
\typeout{** loaded for the language `#1'. Using the pattern for}%
\typeout{** the default language instead.}%
\else
\language=\csname l@#1\endcsname
\fi
#2}}
\providecommand{\BIBdecl}{\relax}
\BIBdecl

\bibitem{wielandallgower2007cbf}
P.~Wieland and F.~Allgöwer, ``Constructive safety using control barrier
  functions,'' \emph{IFAC Proceedings Vols.}, vol.~40, no.~12, pp. 462--467,
  2007.

\bibitem{ames2017control}
A.~D. Ames, X.~Xu, J.~W. Grizzle, and P.~Tabuada, ``Control barrier function
  based quadratic programs for safety critical systems,'' \emph{IEEE Trans. on
  Automatic Control}, vol.~62, no.~8, pp. 3861--3876, 2017.

\bibitem{Xiao2019HOCBF}
W.~Xiao and C.~Belta, ``Control barrier functions for systems with high
  relative degree,'' in \emph{IEEE 58th Conference on Decision and Control
  (CDC)}, 2019, pp. 474--479.

\bibitem{cortez2019cbfmechsys}
W.~S. Cortez, D.~Oetomo, C.~Manzie, and P.~Choong, ``Control barrier functions
  for mechanical systems: Theory and application to robotic grasping,'' in
  \emph{IEEE Trans. on Control Systems Tech.}, 2019, pp. 1--16.

\bibitem{Chen2018Obstacle}
Y.~Chen, H.~Peng, and J.~Grizzle, ``Obstacle avoidance for low-speed autonomous
  vehicles with barrier function,'' \emph{IEEE Transactions on Control Systems
  Technology}, vol.~26, no.~1, pp. 194--206, 2018.

\bibitem{Garg2021Robust}
K.~Garg and D.~Panagou, ``Robust control barrier and control lyapunov functions
  with fixed-time convergence guarantees,'' in \emph{2021 American Control
  Conference (ACC)}, 2021, pp. 2292--2297.

\bibitem{Chen2021Guaranteed}
Y.~Chen, A.~Singletary, and A.~D. Ames, ``Guaranteed obstacle avoidance for
  multi-robot operations with limited actuation: A control barrier function
  approach,'' \emph{IEEE Control Systems Letters}, vol.~5, no.~1, pp. 127--132,
  2021.

\bibitem{Jankovic2021Collision}
M.~Jankovic and M.~Santillo, ``Collision avoidance and liveness of multi-agent
  systems with cbf-based controllers,'' in \emph{2021 60th IEEE Conference on
  Decision and Control (CDC)}, 2021, pp. 6822--6828.

\bibitem{Xu2018Safe}
B.~Xu and K.~Sreenath, ``Safe teleoperation of dynamic uavs through control
  barrier functions,'' in \emph{2018 IEEE International Conference on Robotics
  and Automation (ICRA)}, 2018, pp. 7848--7855.

\bibitem{Khan2020Cascaded}
M.~Khan, M.~Zafar, and A.~Chatterjee, ``Barrier functions in cascaded
  controller: Safe quadrotor control,'' in \emph{2020 American Control
  Conference (ACC)}, 2020, pp. 1737--1742.

\bibitem{Black2022ffcbf}
M.~Black, M.~Jankovic, A.~Sharma, and D.~Panagou, ``Future-focused control
  barrier functions for autonomous vehicle control,'' in \emph{2023 American
  Control Conference (ACC), to appear}, 2023.

\bibitem{Yaghoubi2021RiskBoundedCBF}
S.~Yaghoubi, G.~Fainekos, T.~Yamaguchi, D.~Prokhorov, and B.~Hoxha,
  ``Risk-bounded control with kalman filtering and stochastic barrier
  functions,'' in \emph{2021 60th IEEE Conference on Decision and Control
  (CDC)}, 2021, pp. 5213--5219.

\bibitem{Wang2018Permissive}
L.~Wang, D.~Han, and M.~Egerstedt, ``Permissive barrier certificates for safe
  stabilization using sum-of-squares,'' in \emph{2018 Annual American Control
  Conference (ACC)}, 2018, pp. 585--590.

\bibitem{clark2021verification}
A.~Clark, ``Verification and synthesis of control barrier functions,'' in
  \emph{60th IEEE Conference on Decision and Control (CDC)}, 2021, pp.
  6105--6112.

\bibitem{pond2022fast}
E.~Pond and M.~Hale, ``Fast verification of control barrier functions via
  linear programming,'' \emph{arXiv preprint arXiv:2212.00598}, 2022.

\bibitem{Breeden2021InputConstraints}
J.~Breeden and D.~Panagou, ``High relative degree control barrier functions
  under input constraints,'' in \emph{2021 60th IEEE Conference on Decision and
  Control (CDC)}, 2021, pp. 6119--6124.

\bibitem{Xiao2022Adaptive}
W.~Xiao, C.~Belta, and C.~G. Cassandras, ``Adaptive control barrier
  functions,'' \emph{IEEE Transactions on Automatic Control}, vol.~67, no.~5,
  pp. 2267--2281, 2022.

\bibitem{Garg2022Fixed}
K.~Garg, E.~Arabi, and D.~Panagou, ``Fixed-time control under spatiotemporal
  and input constraints: A quadratic programming based approach,''
  \emph{Automatica}, vol. 141, p. 110314, 2022.

\bibitem{Garg2019Control}
K.~Garg and D.~Panagou, ``Control-lyapunov and control-barrier functions based
  quadratic program for spatio-temporal specifications,'' in \emph{58th IEEE
  Conf. on Decision and Control (CDC)}, 2019, pp. 1422--1429.

\bibitem{Black2020Quadratic}
M.~Black, K.~Garg, and D.~Panagou, ``A quadratic program based control
  synthesis under spatiotemporal constraints and non-vanishing disturbances,''
  in \emph{2020 59th IEEE Conference on Decision and Control (CDC)}, 2020, pp.
  2726--2731.

\bibitem{Shao2021Tracking}
K.~Shao, J.~Zheng, H.~Wang, X.~Wang, R.~Lu, and Z.~Man, ``Tracking control of a
  linear motor positioner based on barrier function adaptive sliding mode,''
  \emph{IEEE Transactions on Industrial Informatics}, vol.~17, no.~11, pp.
  7479--7488, 2021.

\bibitem{Lindemann2019CBFSTL}
L.~Lindemann and D.~V. Dimarogonas, ``Control barrier functions for signal
  temporal logic tasks,'' \emph{IEEE Control Systems Letters}, vol.~3, no.~1,
  pp. 96--101, 2019.

\bibitem{Yang2020STL}
G.~Yang, C.~Belta, and R.~Tron, ``Continuous-time signal temporal logic
  planning with control barrier functions,'' in \emph{2020 American Control
  Conference (ACC)}, 2020, pp. 4612--4618.

\bibitem{Srinivasan2021Control}
M.~Srinivasan and S.~Coogan, ``Control of mobile robots using barrier functions
  under temporal logic specifications,'' \emph{IEEE Transactions on Robotics},
  vol.~37, no.~2, pp. 363--374, 2021.

\bibitem{Cortez2022RobustMultiple}
W.~S. Cortez, X.~Tan, and D.~V. Dimarogonas, ``A robust, multiple control
  barrier function framework for input constrained systems,'' \emph{IEEE
  Control Systems Letters}, vol.~6, pp. 1742--1747, 2022.

\bibitem{Glotfelter2017Nonsmooth}
P.~Glotfelter, J.~Cortes, and M.~Egerstedt, ``Nonsmooth barrier functions with
  applications to multi-robot systems,'' \emph{IEEE Control Systems Letters},
  vol.~1, pp. 310--315, 10 2017.

\bibitem{huang2020switched}
Y.~Huang and Y.~Chen, ``Switched control barrier function with applications to
  vehicle safety control,'' in \emph{Dynamic Systems and Control Conference},
  vol. 84270.\hskip 1em plus 0.5em minus 0.4em\relax American Society of
  Mechanical Engineers, 2020, p. V001T15A002.

\bibitem{Machida2021ConsensusCBF}
M.~Machida and M.~Ichien, ``Consensus-based control barrier function for
  swarm,'' in \emph{2021 IEEE International Conference on Robotics and
  Automation (ICRA)}, 2021, pp. 8623--8628.

\bibitem{cai2021safe}
Z.~Cai, H.~Cao, W.~Lu, L.~Zhang, and H.~Xiong, ``Safe multi-agent reinforcement
  learning through decentralized multiple control barrier functions,''
  \emph{arXiv preprint arXiv:2103.12553}, 2021.

\bibitem{Ma2021Generalized}
H.~Ma, J.~Chen, S.~Eben, Z.~Lin, Y.~Guan, Y.~Ren, and S.~Zheng, ``Model-based
  constrained reinforcement learning using generalized control barrier
  function,'' in \emph{2021 IEEE/RSJ International Conference on Intelligent
  Robots and Systems (IROS)}, 2021, pp. 4552--4559.

\bibitem{breeden2022compositions}
J.~Breeden and D.~Panagou, ``Compositions of multiple control barrier functions
  under input constraints,'' \emph{arXiv preprint arXiv:2210.01354}, 2022.

\bibitem{Black2022Fixed}
M.~Black, E.~Arabi, and D.~Panagou, ``Fixed-time parameter adaptation for safe
  control synthesis,'' \emph{arXiv preprint arXiv:2204.10453}, 2022.

\bibitem{Lopez2021RaCBF}
B.~T. Lopez, J.-J.~E. Slotine, and J.~P. How, ``Robust adaptive control barrier
  functions: An adaptive and data-driven approach to safety,'' \emph{IEEE
  Control Systems Letters}, vol.~5, no.~3, pp. 1031--1036, 2021.

\bibitem{Taylor2020Adaptive}
A.~J. Taylor and A.~D. Ames, ``Adaptive safety with control barrier
  functions,'' in \emph{2020 American Control Conference (ACC)}, 2020, pp.
  1399--1405.

\bibitem{black2022adaptation}
M.~Black and D.~Panagou, ``Adaptation for validation of a consolidated control
  barrier function based control synthesis,'' \emph{arXiv preprint
  arXiv:2209.08170}, 2022.

\bibitem{bazaraa2013nonlinear}
M.~S. Bazaraa, H.~D. Sherali, and C.~M. Shetty, \emph{Nonlinear programming:
  theory and algorithms}.\hskip 1em plus 0.5em minus 0.4em\relax John Wiley \&
  Sons, 2013.

\bibitem{bertsekas1979convexification}
D.~P. Bertsekas, ``Convexification procedures and decomposition methods for
  nonconvex optimization problems,'' \emph{Journal of Optimization Theory and
  Applications}, vol.~29, no.~2, pp. 169--197, 1979.

\bibitem{Fazlyab2018Prediction}
M.~Fazlyab, S.~Paternain, V.~M. Preciado, and A.~Ribeiro,
  ``Prediction-correction interior-point method for time-varying convex
  optimization,'' \emph{IEEE Trans. on Automatic Control}, vol.~63, no.~7, pp.
  1973--1986, 2018.

\bibitem{breeden2021robust}
J.~Breeden and D.~Panagou, ``Robust control barrier functions under high
  relative degree and input constraints for satellite trajectories,''
  \emph{Automatica}, (to be published) 2023.

\bibitem{nguyen2016exponential}
Q.~Nguyen and K.~Sreenath, ``Exponential control barrier functions for
  enforcing high relative-degree safety-critical constraints,'' in \emph{2016
  American Control Conference (ACC)}.\hskip 1em plus 0.5em minus 0.4em\relax
  IEEE, 2016, pp. 322--328.

\bibitem{Jankovic2022Multi}
M.~Jankovic, M.~Santillo, and Y.~Wang, ``Multi-agent systems with cbf-based
  controllers -- collision avoidance and liveness from instability,''
  \emph{arXiv preprint arXiv:2207.04915}, 2022.

\bibitem{xiao2022sufficient}
W.~Xiao, C.~A. Belta, and C.~G. Cassandras, ``Sufficient conditions for
  feasibility of optimal control problems using control barrier functions,''
  \emph{Automatica}, vol. 135, p. 109960, 2022.

\bibitem{wang2022suboptimal}
S.~Wang, K.~Shi, T.~Huang \emph{et~al.}, ``Suboptimal safety-critical control
  for continuous systems using prediction-correction online optimization,''
  \emph{arXiv preprint arXiv:2203.15305}, 2022.

\bibitem{khalil2002nonlinear}
H.~K. Khalil, \emph{Nonlinear systems}.\hskip 1em plus 0.5em minus 0.4em\relax
  Prentice hall Upper Saddle River, NJ, 2002, vol.~3.

\bibitem{Rajamani2012VDC}
R.~Rajamani, \emph{Vehicle Dynamics and Control}.\hskip 1em plus 0.5em minus
  0.4em\relax Springer US, 2012.

\bibitem{Xiao2022HighOrder}
W.~Xiao and C.~Belta, ``High-order control barrier functions,'' \emph{IEEE
  Transactions on Automatic Control}, vol.~67, no.~7, pp. 3655--3662, 2022.

\end{thebibliography}

\begin{IEEEbiography}[{\includegraphics[width=1in,height=1.25in,clip,keepaspectratio]{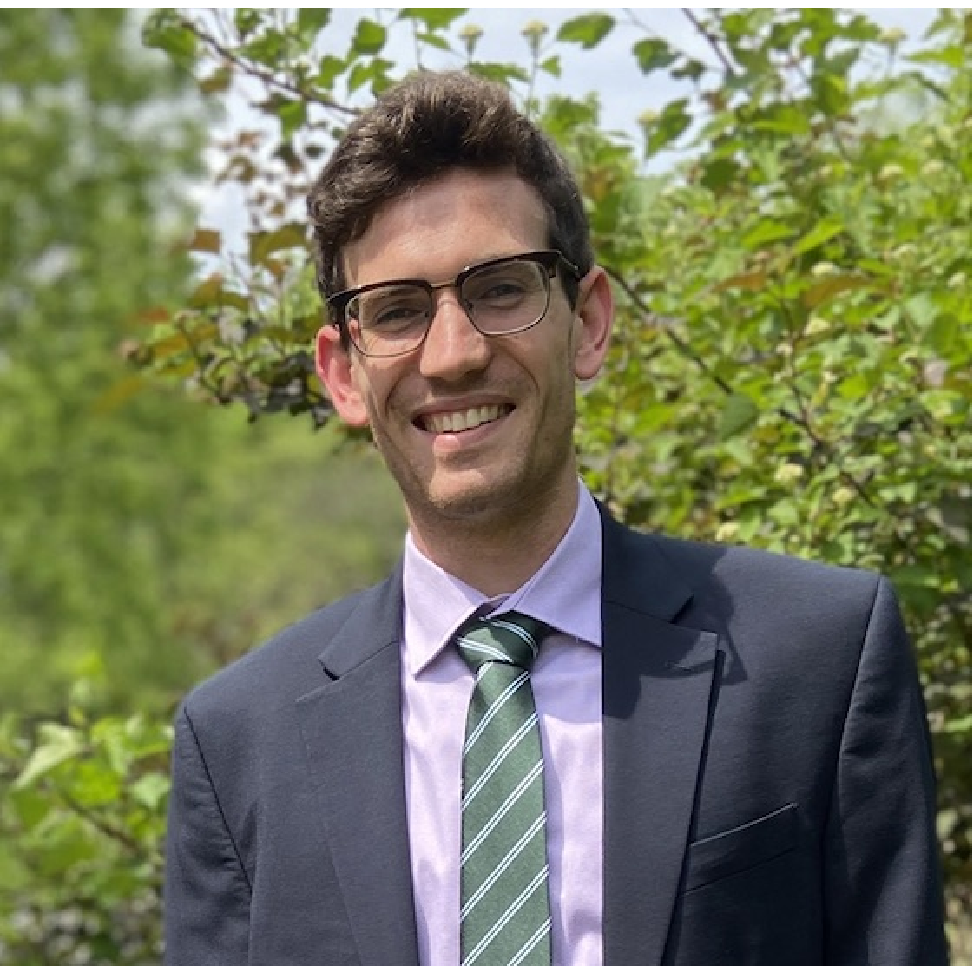}}]{Mitchell Black} (M'20) grew up in Brunswick, Maine, USA. He received his B.S. from Tufts University in 2016 in both Mechanical Engineering and Astrophysics, and his M.S.E. in Aerospace Engineering from the University of Michigan in 2017, where he has been a Ph.D. student in Professor Dimitra Panagou's Distributed Aerospace Systems and Control Lab since 2020. In the interim (2017-2020), he worked as a Research Scientist at Michigan Aerospace Corporation in Ann Arbor, MI.

Mitchell was the 2016 recipient of the National Collegiate Athletic Association's Walter Byers Postgraduate Scholarship for combined achievement in the classroom and as a member of the Tufts University Track \& Field team. He has continued to run throughout graduate school and in 2021 competed at the US Olympic Track \& Field Trials in the 800m, where he placed 23rd.
\end{IEEEbiography}

\begin{IEEEbiography}[{\includegraphics[width=1in,height=1.25in,clip,keepaspectratio]{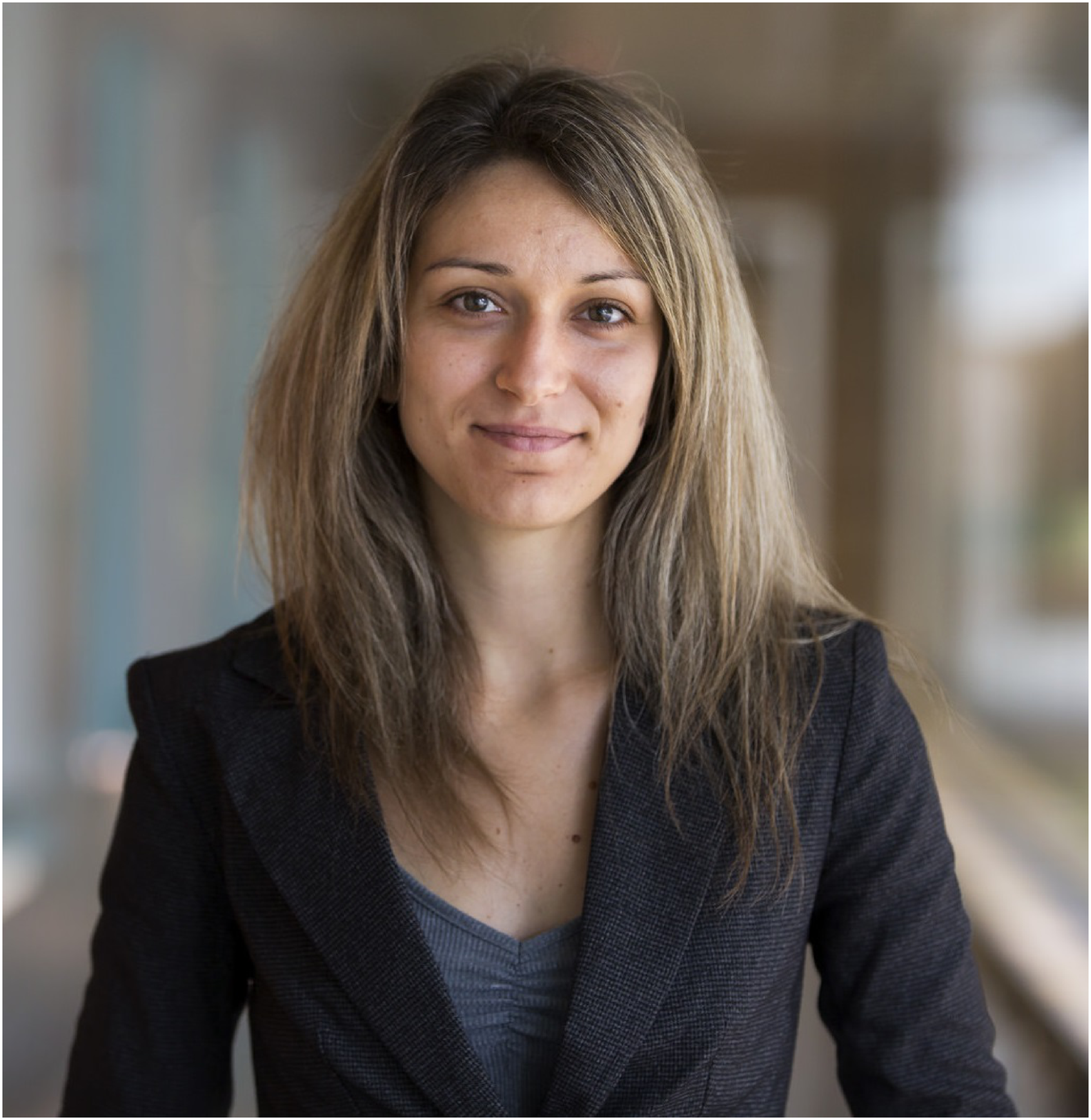}}]{Dimitra Panagou} received the Diploma and PhD degrees in Mechanical Engineering from the National Technical University of Athens, Greece, in 2006 and 2012, respectively. In September 2014 she joined the Department of Aerospace Engineering, University of Michigan as an Assistant Professor. Since July 2022 she is an Associate Professor with the newly established Department of Robotics, with a courtesy “dry” appointment with the Department of Aerospace Engineering, University of Michigan. Prior to joining the University of Michigan, she was a postdoctoral research associate with the Coordinated Science Laboratory, University of Illinois, Urbana-Champaign (2012-2014), a visiting research scholar with the GRASP Lab, University of Pennsylvania (June 2013, Fall 2010) and a visiting research scholar with the University of Delaware, Mechanical Engineering Department (Spring 2009). Her research program spans the areas of nonlinear systems and control; multi-agent systems and networks; motion and path planning; human-robot interaction; navigation, guidance, and control of aerospace vehicles. She is particularly interested in the development of provably-correct methods for the safe and secure (resilient) operation of autonomous systems in complex missions, with applications in robot/sensor networks and multi-vehicle systems (ground, marine, aerial, space). She is a recipient of the NASA Early Career Faculty Award, the AFOSR Young Investigator Award, the NSF CAREER Award, and a Senior Member of the IEEE and the AIAA.

\end{IEEEbiography}

% \begin{IEEEbiography}[{\includegraphics[width=1in,height=1.25in,clip,keepaspectratio]{a3.png}}]{Third C. Author, Jr.} (M'87) received the B.S. degree in mechanical 
% engineering from National Chung Cheng University, Chiayi, Taiwan, in 2004 
% and the M.S. degree in mechanical engineering from National Tsing Hua 
% University, Hsinchu, Taiwan, in 2006. He is currently pursuing the Ph.D. 
% degree in mechanical engineering at Texas A{\&}M University, College 
% Station, TX, USA.

% From 2008 to 2009, he was a Research Assistant with the Institute of 
% Physics, Academia Sinica, Tapei, Taiwan. His research interest includes the 
% development of surface processing and biological/medical treatment 
% techniques using nonthermal atmospheric pressure plasmas, fundamental study 
% of plasma sources, and fabrication of micro- or nanostructured surfaces. 

% Mr. Author's awards and honors include the Frew Fellowship (Australian 
% Academy of Science), the I. I. Rabi Prize (APS), the European Frequency and 
% Time Forum Award, the Carl Zeiss Research Award, the William F. Meggers 
% Award and the Adolph Lomb Medal (OSA).
% \end{IEEEbiography}

\end{document}